\documentclass[11pt,a4paper]{article}
\usepackage{verbatim}
\usepackage{stix}

\usepackage[T1]{fontenc}
\usepackage[a4paper,includeheadfoot,margin=2.54cm]{geometry}
\usepackage{amsthm, amsmath}
\usepackage{xfrac}

\theoremstyle{plain} 
\newtheorem{theorem}{Theorem}
\newtheorem{corollary}{Corollary}
\newtheorem{lemma}{Lemma}
\newtheorem{proposition}{Proposition}

\theoremstyle{definition}
\newtheorem{definition}{Definition}
\newtheorem{remark}{Remark}

\newtheorem{eexample}{Example}

\usepackage{algorithm}
\usepackage{algorithmic}
\newcommand{\algAnd}{\textbf{and} }
\newcommand{\algOr}{\textbf{or} }

\renewcommand{\algorithmicrequire}{\textbf{Input:}}
\renewcommand{\algorithmicensure}{\textbf{Output:}}
\makeatletter
\newenvironment{breakablealgorithm}
  {
   \begin{center}
     \refstepcounter{algorithm}
     \hrule height.8pt depth0pt \kern2pt
     \renewcommand{\caption}[2][\relax]{
       {\raggedright\textbf{\fname@algorithm~\thealgorithm} ##2\par}%
       \ifx\relax##1\relax 
         \addcontentsline{loa}{algorithm}{\protect\numberline{\thealgorithm}##2}%
       \else 
         \addcontentsline{loa}{algorithm}{\protect\numberline{\thealgorithm}##1}%
       \fi
       \kern2pt\hrule\kern2pt
     }
  }{
     \kern2pt\hrule\relax
   \end{center}
  }
\makeatother

\usepackage{amsfonts}
\usepackage[shortlabels]{enumitem}

\newcommand{\makeblue}[1]{{\color{blue}#1}}
\newcommand{\new}[1]{#1}

\newcommand{\C}{\mathbb{C}}
\newcommand{\N}{\mathbb{N}}
\newcommand{\R}{\mathbb{R}}
\newcommand{\T}{\mathbb{T}}

\newcommand{\Z}{\mathbb{Z}}
\newcommand{\Rnn}{\R^{n\times n}}
\newcommand{\Cnn}{\C^{n\times n}}
\newcommand{\CNN}[1]{\C^{{#1}\times {#1}}}
\newcommand{\Cn}{\C^n}

\newcommand{\bigO}{{\mathcal{O}}}
\renewcommand{\d}[1][x]{\,\operatorname{d}\!#1}

\newcommand{\ddt}{\frac{\d[]}{\d[t]}}

\newcommand{\eps}{\varepsilon}

\newcommand{\mA}{\mathbf{A}}
\newcommand{\mB}{\mathbf{B}}
\newcommand{\mC}{\mathbf{C}}
\newcommand{\mD}{\mathbf{D}}
\newcommand{\mE}{\mathbf{E}}
\newcommand{\mF}{\mathbf{F}}
\newcommand{\mG}{\mathbf{G}}

\newcommand{\mI}{\mathbf{I}}
\newcommand{\mJ}{\mathbf{J}}

\newcommand{\mL}{\mathbf{L}}
\newcommand{\mM}{\mathbf{M}}
\newcommand{\mN}{\mathbf{N}}

\newcommand{\mP}{\mathbf{P}}
\newcommand{\mQ}{\mathbf{Q}}

\newcommand{\mR}{\mathbf{R}}
\newcommand{\mS}{\mathbf{S}}
\newcommand{\mT}{\mathbf{T}}
\newcommand{\mU}{\mathbf{U}}
\newcommand{\mV}{\mathbf{V}}
\newcommand{\mW}{\mathbf{W}}
\newcommand{\mX}{\mathbf{X}}
\newcommand{\mY}{\mathbf{Y}}
\newcommand{\mZ}{\mathbf{Z}}

\newcommand{\msqrtE}{\mE^{1/2}}
\newcommand{\msqrtEinv}{(\msqrtE)^{-1}}

\newcommand{\tA}{\widetilde \mA}
\newcommand{\tB}{\widetilde \mB}
\newcommand{\tD}{\widetilde \mD}
\newcommand{\tE}{\widetilde \mE}
\newcommand{\tG}{\widetilde \mG}
\newcommand{\tJ}{\widetilde \mJ}
\newcommand{\tN}{\widetilde \mN}
\newcommand{\tR}{\widetilde \mR}
\newcommand{\tSigma}{\widetilde \Sigma}
\newcommand{\tX}{\widetilde \mX}
\newcommand{\tY}{\widetilde \mY}

\newcommand{\tn}{\tilde n}

\newcommand{\mAH}{\mA_H}
\newcommand{\mAS}{\mA_S}

\newcommand{\mCH}{\mC_H}
\newcommand{\mCS}{\mC_S}

\newcommand{\mHC}{m_{HC}}
\newcommand{\xk}{{w_k}}

\newcommand{\bbH}{\mathbb{H}}
\newcommand{\Hn}{\mathbb{H}_n} 
\newcommand{\PDHn}{\mathbb{H}_n^{>}} 
\newcommand{\PSDHn}{\mathbb{H}_n^{\geq}} 

\newcommand{\ip}[2]{\langle {#1}\ ,\, {#2} \rangle}

\DeclareMathOperator{\diver}{div}
\DeclareMathOperator{\diag}{diag}
\DeclareMathOperator{\kernel}{ker}

\DeclareMathOperator{\rank}{rank}

\newcommand{\cE}{\mathcal{E}}
\newcommand{\cG}{\mathcal{G}}
\newcommand{\cH}{\mathcal{H}}

\newcommand{\cV}{\mathcal{V}}

\newcommand{\bx}{x}

\begin{document}

\title{Hypocoercivity and controllability in linear semi-dissipative \new{Hamiltonian} ODEs and DAEs}

\author{Franz Achleitner\thanks{Vienna University of Technology, Institute of Analysis and Scientific Computing, Wiedner Hauptstr. 8-10, A-1040 Wien, Austria, franz.achleitner@tuwien.ac.at},
Anton Arnold\thanks{Vienna University of Technology, Institute of Analysis and Scientific Computing, Wiedner Hauptstr. 8-10, A-1040 Wien, Austria, anton.arnold@tuwien.ac.at},
and Volker Mehrmann\thanks{Technische Universit\"at Berlin, Institut f.~Mathematik,  MA 4-5, Stra\ss{}e des 17.~Juni 136, D-10623 Berlin, mehrmann@math.tu-berlin.de}
} 


\maketitle

\begin{abstract}
For the classes of finite dimensional linear \new{time-invariant} semi-dissipative Hamiltonian ordinary differential equations \new{and} differential-algebraic equations, stability and hypocoercivity are discussed and related to concepts from control theory.
On the basis of staircase forms the solution behavior is characterized and
connected to the hypocoercivity index of these evolution equations.
The results are applied to two infinite dimensional flow problems.
\end{abstract}

\section{Introduction}
\label{sec:intro}
A detailed analysis of the stability of dynamical systems of evolution equations (finite or infinite-di\-men\-sion\-al) is still very problem dependent and computationally challenging, see~\cite{Adr95,DieRV97,DieV02b,HiPr10}.
In view of these challenges it is important to use structural information of the dynamical system to characterize stability, asymptotic stability and the transient behavior.
In this paper we consider these questions for two classes of finite-dimensional linear systems, although we have in mind to extend these results to the infinite-dimensional case and will do so for some examples.

The first class are ordinary differential equations (ODEs) of the form
\begin{equation}\label{ODE:A}
 \dot{x}(t) = \mA x(t)\,, \qquad t\geq 0\,,
\end{equation}
for some function $x:[0,\infty)\to\C^n$ and a constant matrix~$\mA\in\Cnn$.
The second, and more general class are differential-algebraic equations (DAEs) of the form
\begin{equation}\label{DAE:EA}
\mE \dot{x}(t) = \mA x(t)\,, \qquad t\geq 0\,,
\end{equation}
for a pair~$(\mE,\mA)$ of constant matrices~$\mE,\mA\in\Cnn$ with $\mE=\mE^H$ positive semi-definite.
Note that if~$\mE$ is positive definite, then it has a matrix square root~$\msqrtE$.
Then the DAE~\eqref{DAE:EA} takes the form~\eqref{ODE:A} by changing the variables~$y:= \msqrtE x$ and scaling the equation by $\msqrtEinv$, such that
\begin{equation}\label{ODE:A-tilde}
\dot{y}(t) =\tA y
\qquad\text{where }
\tA :=\msqrtEinv\mA \msqrtEinv \ .
\end{equation}
However, if $\mE$ is singular then the behavior of the two systems is fundamentally different.

We use the following notation:
The conjugate transpose (transpose) of a matrix $\mC$ is denoted by $\mC^H$ ($\mC^\top$).
The set of Hermitian matrices in $\Cnn$ is denoted by $\bbH_n$.
Positive definiteness (semi-definiteness) of $\mC\in\bbH_n$ is denoted by $\mC>0$ ($\mC\geq 0$).
The set of all positive definite (semi-definite) Hermitian matrices in $\Cnn$ is denoted by $\bbH_n^>$ ($\bbH_n^\geq$).

Writing $\mA$ as the sum of its Hermitian part $\mAH =\tfrac12 (\mA +\mA^H)$ and skew-Hermitian part $\mAS =\tfrac12 (\mA -\mA^H)$, we have the following definition.
\begin{definition}[{Definition 4.1.1 of~\cite{Be18}}]
A matrix~$\mA\in\Cnn$  is called~\emph{dissipative} (resp. \emph{semi-dissipative}) if the Hermitian part~$\mAH$ is negative definite (resp. negative semi-definite).

For a (semi-)dissipative matrix~$\mA\in\Cnn$, the associated ODE~\eqref{ODE:A} is called \emph{(semi-)dissipative Hamiltonian ODE}.

For a (semi-)dissipative matrix~$\mA\in\Cnn$ and positive semi-definite Hermitian matrix~$\mE\in\Cnn$, the associated DAE~\eqref{DAE:EA} is called \emph{(semi-)dissipative Hamiltonian DAE}.
\end{definition}
Since $\mE=\mE^H\geq 0$, note that the above definition for the DAE~\eqref{DAE:EA} is independent of~$\mE$ (see Remark~\ref{remark:congruence} and Theorem~\ref{thm:singind}~\ref{lambda:non-pos} below).
The notion (semi-)dissipative Hamiltonian is motivated by the fact that if $\mAH=0$, then~\eqref{ODE:A} is a Hamiltonian system with Hamiltonian $\mathsf H=1/2\ x^H x$, see~\cite{MehMW18}.
In the following (to avoid too many indices), we write semi-dissipative matrices~$\mA$ in the form~$\mA =\mJ -\mR$ with a skew-Hermitian matrix~$\mJ=\mAS$ and a positive semi-definite Hermitian matrix $\mR =-\mAH$.

\medskip
The notion \emph{hypocoercivity} was introduced 
in \cite{Vi09} for equations (mostly partial differential equations) of the form $\ddt{x}=-\mC x$ on some Hilbert space~$\cH$, where the linear generator~$\mC$ is not coercive, but where solutions still exhibit exponential decay in time.
More precisely, for hypocoercive generators~$\mC$ there exist constants $\lambda>0$ and $c\ge1$, such that
\begin{equation}\label{exp-decay}
  \|e^{-\mC t} x_0\|_{\widetilde \cH}
  \le c\,e^{-\lambda t} \|x_0\|_{\widetilde \cH}
  \qquad
  \mbox{\rm for all}\,x_0\in \widetilde \cH\,,
  \qquad
  t \geq 0\,,
\end{equation}
where $\widetilde \cH$ is another Hilbert space, densely embedded in $(\ker \mC)^\perp\subset \cH$.
Often, this phenomenon is also called hypocoercivity.

The long-time behavior of many systems exhibiting hypocoercivity has been studied in recent years, including Fokker--Planck equations~\cite{AASt15, ArEr14, Vi09}, kinetic equations~\cite{DoMoSc09,DoMoSc15}, and reaction-transport equations of BGK-type~\cite{AAC16, AAC18}.
Determining the sharp (i.e. maximal) exponential decay rate~$\lambda$ was an issue in some of these works, in particular \cite{AAC16, AAC18, ArEr14}.
But finding at the same time the smallest multiplicative constant~$c\ge1$, is a rather recent topic, e.g. see~\cite{AAS19,AS21,GuMo16,GuMo16E}.

\new{
The \emph{hypocoercivity index (HC-index)} is a measure for the structural complexity of the intertwining between the Hermitian and skew-Hermitian part of hypocoercive generators, and it determines the short time behavior of its solutions~\cite{AAC20,ArEr14}.
}

\smallskip
\new{
In this article we shall present the following main results:
\begin{itemize}
\item 
We show that the HC-index of the matrix $\mA$ from~\eqref{ODE:A} can be read directly from a \emph{staircase form} of $\mA$. 
Since the latter is obtained via a unitary transformation, this approach allows to compute the HC-index $\mHC$ with a numerically stable algorithm.
\item 
For particular control systems, we relate the HC-index to the \emph{controllability index} $m_C$ as $\mHC = m_C-1$.
\item 
For semi-dissipative Hamiltonian ODEs~\eqref{ODE:A}, we give an iterative construction of a strict Lyapunov functional (in form of a weighted Euclidean norm), where the number of iteration steps matches $\mHC$.
\item 
For semi-dissipative Hamiltonian DAEs~\eqref{DAE:EA}, we give a unitary transformation to a staircase form and an equivalence transformation to an \emph{almost Kronecker canonical form}. 
This is our basis to define (apparently for the first time) the HC-index for linear DAEs. 
It matches the HC-index for the implicit ODE system that describes the (semi-)dissipative dynamics of~\eqref{DAE:EA}, if consistent initial conditions are prescribed.
\item 
Extending the ODE-case, this HC-index also characterizes the short time behavior of semi-dissipative Hamiltonian DAEs (in the $\mE$-weighted norm, with the matrix $\mE$ given in~\eqref{DAE:EA}).
\end{itemize}
}

\medskip
In Section~\ref{sec:ODE}, we discuss linear systems of ordinary differential equations (ODEs) and recall the staircase algorithm for pairs of matrices where at least one is skew-Hermitian.
Then, we introduce the notion of (hypo)coercive matrices, define the HC-index and relate it to the \emph{controllability index} of control theory.
The hypocoercivity index (resp. controllability index) appears in the construction of strict Lyapunov functionals, see e.g.~\cite{AAC20,Mi09thesis}.

In Section~\ref{sec:DAE} the results from Section~\ref{sec:ODE} are extended to DAEs of the form~\eqref{DAE:EA}, where in particular the extension of the hypocoercivity as well as the Lyapunov theory is discussed on the basis of staircase forms.
In Section~\ref{sec:example} the results are applied to Stokes flow and acoustic waves in pipeline networks.


\section{Linear time-invariant ODE systems}
\label{sec:ODE}

\subsection{Canonical and condensed forms for semi-dissipative \new{Hamiltonian} ODEs}
\label{ssec:ODE:CF}
An equivalence relation on a set allows to define a canonical form, or a condensed form under unitary transformations. 
We use the following four equivalence relations in matrix analysis, see e.g.~\S0.11 in~\cite{HoJo13} or \S4.7 in~\cite{Be18}:
Two matrices~$\mA, \tA\in\Cnn$ are called
\begin{itemize}
\item
\emph{equivalent} if there exist nonsingular matrices~$\mL,\mZ\in\Cnn$ such that $\tA=\mL\mA\mZ$;
\item
\emph{similar} if there exists an invertible matrix~$\mS\in\Cnn$ such that $\tA=\mS\mA\mS^{-1}$;
\item
\emph{congruent} (or \emph{*congruent}) if there exists a nonsingular matrix~$\mQ\in\Cnn$ such that $\tA=\mQ\mA\mQ^H$;
\item
\emph{unitarily congruent} (or \emph{unitarily *congruent} or \emph{unitarily similar}) if there exists a unitary matrix~$\mP\in\Cnn$ such that $\tA=\mP\mA\mP^H$.
We prefer the notion~\emph{unitary congruence} since congruence preserves e.g.~the inertia of the Hermitian part of matrices.
\end{itemize}
Each of these relations is an \emph{equivalence relation}, denoted as $\mA\sim \tA$, and allows to define a canonical form or condensed form.
\begin{remark}\label{remark:congruence}
Consider a semi-dissipative Hamiltonian DAE~\eqref{DAE:EA} with positive definite Hermitian matrix~$\mE$.
Using an equivalence transformation yields an ODE~\eqref{ODE:A} with modified matrix~$\mA$, which may or may not be semi-dissipative.
For example, 
scaling Equation~\eqref{DAE:EA} by~$\mE^{-1}$ yields an ODE of the form:
\begin{equation}\label{ODE:E-1A}
\dot{x}(t) =\mE^{-1}\mA x(t)\ .
\end{equation}
But using the congruence transformation of~$\mA$ as in~\eqref{ODE:A-tilde} yields an ODE~\eqref{ODE:A-tilde} which is again semi-dissipative.
By contrast, if~$\mE$ is singular then the behavior of DAE systems is fundamentally different (compared with ODE systems).
\end{remark}
To illustrate the difference between the equivalence transformation $(\mE,\mA)\mapsto(\mI,\mE^{-1}\mA)$ and the congruence transformation $(\mE,\mA)\mapsto(\mI,\msqrtEinv\mA \msqrtEinv)$, let us consider the following example:
\begin{eexample} \label{ex:comparison:equiv+cong}
Consider the semi-dissipative Hamiltonian DAE~\eqref{DAE:EA} with matrices
\begin{align}\label{DAE:EA:example2x2}
\mA
&:=\begin{bmatrix} -1 & 1 \\ -1 & 0 \end{bmatrix} ,
 & %
\mE
&:=\frac{16}9
\begin{bmatrix} \sfrac54 & -1 \\ -1 & \sfrac54 \end{bmatrix} ,
 & %
\mE^{-1}
&=\begin{bmatrix} \sfrac54 & 1 \\ 1 & \sfrac54 \end{bmatrix} ,
 & %
\msqrtEinv
&=\begin{bmatrix} 1 & \sfrac12 \\ \sfrac12 & 1 \end{bmatrix} .
\end{align}
The semi-dissipative matrix~$\mA$ is stable and has eigenvalues~$\lambda_\pm =(-1\pm i\sqrt{3})/2$.
Its Hermitian part~$\mAH =\diag({-1}, 0)$ is only semi-definite.
Then, we compute the matrices $\mE^{-1}\mA$ and $\tA:=\msqrtEinv\mA\msqrtEinv$ and their Hermitian parts as
\begin{align*}
\mE^{-1}\mA
&=\begin{bmatrix} -\sfrac94 & \sfrac54 \\ -\sfrac94 & 1 \end{bmatrix} ,
 & %
\big(\mE^{-1}\mA\big)_H
&=\begin{bmatrix} -\sfrac94 & -\sfrac12 \\ -\sfrac12 & 1 \end{bmatrix} ,
 & %
\tA
&=\begin{bmatrix} -1 & \sfrac14 \\ -\sfrac54 & -\sfrac14 \end{bmatrix} ,
 & %
(\tA)_H
&=\begin{bmatrix} -1 & -\sfrac12 \\ -\sfrac12 & -\sfrac14 \end{bmatrix} .
\end{align*}
The matrices~$\mE^{-1}\mA$ and~$\tA$ are similar, since $\mE^{-1}\mA =\msqrtEinv\tA\msqrtE$; hence, they share the same eigenvalues~$(-5\pm i\sqrt{11})/8$ and are stable.
However, their Hermitian parts are very different. 
In fact, $(\mE^{-1}\mA)_H$ is indefinite, whereas $(\tA)_H$ inherits the negative semi-definiteness from $\mAH$ due to the congruence transformation.
\end{eexample}
%
%
In the following we make frequent use of the condensed (staircase) form of matrix pairs, see~\cite{VD79}, adapted to the special structure of~$\mJ$ and~$\mR$.
\begin{lemma}[Staircase form for $(\mJ,\mR)$] \label{lem:SF}
Let $\mJ\in\Cnn$ be a skew-Hermitian matrix, and $\mR\in\Cnn$ be a nonzero Hermitian matrix.
Then there exists a unitary matrix $\mP\in\Cnn$, such that~$\mP \mJ \mP^H$ and~$\mP \mR \mP^H$ are block tridiagonal matrices of the form
\begin{equation} \label{matrices:staircase:J-R}
\begin{split}
\mP\ \mJ\ \mP^H
&= \begin{array}{l}
\left[ \begin{array}{ccccccc|c}
 \mJ_{1,1} & -\mJ_{2,1}^H & & & \cdots & & 0 & 0\\
 \mJ_{2,1} & \mJ_{2,2} & -\mJ_{3,2}^H & & & & &\\
  & \ddots & \ddots & \ddots & & & \vdots & \\
  & & \mJ_{k,k-1} & \mJ_{k,k} & -\mJ_{k+1,k}^H & & & \vdots \\
 \vdots & & & \ddots & \ddots & \ddots & & \\
  & & & & \mJ_{s-2,s-3} & \mJ_{s-2,s-2} & -\mJ_{s-1,s-2}^H & \\
 0 & & \cdots & & & \mJ_{s-1,s-2} & \mJ_{s-1,s-1} & 0\\ \hline
 0 & & & \cdots & & & 0 & \mJ_{ss}
\end{array}\right]
 \begin{array}{c}
  n_1\\ n_2\\ \vdots\\ n_k\\ \vdots \\ n_{s-2} \\n_{s-1}\\ n_s
 \end{array} \\
 \quad n_1 \hspace{120pt} n_{s-2} \hspace{20pt} n_{s-1} \hspace{18pt} n_s
\end{array},
\\
\mP\ \mR\ \mP^H
&= \begin{array}{l} \left[ \begin{array}{cc}
     \mR_1 & 0\\
     0 & 0 \\
     \vdots & \vdots\\
     \vdots & \vdots \\
     0 & 0 \end{array}\right]
 \begin{array}{c}
     n_1 \\ n_2 \\ \vdots\\ \vdots \\ n_s \end{array} \\
  \;\; n_1  \;\; n-n_1
\end{array},
\end{split}
\end{equation}
where $n_1 \geq n_2 \geq \cdots \geq n_{s-1} \geq n_s \geq 0$, $n_{s-1}>0$, and $\mR_1\in\C^{n_1,n_1}$ is nonsingular.

If $\mR$ is nonsingular, then~$s=2$ and~$n_2=0$.
For example, $\mP=\mI$, $\mJ_{1,1}=\mJ$ and~$\mR_1=\mR$ is an admissible choice.

If $\mR$ is singular, then $s\geq 3$ and the matrices~$\mJ_{i,i-1}$, $i=2,\ldots,s-1$, in the subdiagonal have full row rank and are of the form
\[
\mJ_{i,i-1} = \begin{bmatrix} \Sigma_{i,i-1} & 0 \end{bmatrix}, \quad
i =2,\ldots, s-1,
\]
%
 %
with nonsingular matrices~$\Sigma_{i,i-1}\in\C^{n_i,n_i}$, moreover $\Sigma_{s-1,s-2}$ is a real-valued diagonal matrix.
\end{lemma}
\begin{proof}
If~$\mR$ is nonsingular, then~$n_1=n$ and we have to choose~$s=2$ and~$n_2=0$ to fit~$\mJ$ into the proposed structure in~\eqref{matrices:staircase:J-R}.

If~$\mR$ is singular, then we give a constructive proof via the following Algorithm~\ref{algorithm:staircase:J-R}.
This algorithm is an adaption of the standard staircase algorithm (to the skew-Hermitian/Hermitian structure of~$(\mJ,\mR)$), see also \cite{BruM07,ByeMX07}.
\begin{breakablealgorithm}
\caption{Staircase Algorithm for pair~$(\mJ,\mR)$}
\label{algorithm:staircase:J-R}
\begin{algorithmic}[1]
\REQUIRE $(\mJ,\mR)$
\newline ----------- \textit{Step 0} -----------
\STATE Perform a (spectral) decomposition of $\mR$ such that
\[
\mR
=\mU_R \begin{bmatrix} \tR_1 & 0 \\ 0 & 0 \end{bmatrix} \mU_R^H ,
\]
with $\mU_R\in \Cnn$ unitary, $\tR_1\in \mathbb C^{n_1, n_1}$ nonsingular.
\STATE Set $\mP := \mU^H_R$, $\tR:=\mU^H_R\ \mR\ \mU_R$,
\[
 \tJ
 := \mU^H_R\ \mJ\ \mU_R
 =:\begin{bmatrix}
   \tJ_{1,1} & \tJ_{1,2} \\
   \tJ_{2,1} & \tJ_{2,2}
   \end{bmatrix},
\]
with $\tJ_{1,1}\in \mathbb C^{n_1, n_1}, \tJ_{2,2}$ skew-Hermitian, and $\tJ_{1,2}=-\tJ_{2,1}^H$.
\newline ----------- \textit{Step 1} -----------
\STATE Perform a singular value decomposition (SVD) of~$\tJ_{2,1}$ such that
\[
 \tJ_{2,1}
 = \mU_{2,1} \begin{bmatrix} \tSigma_{2,1} & 0\\ 0 & 0 \end{bmatrix} \mV^H_{2,1},
\]
with unitary matrices~$\mU_{2,1}$ and~$\mV_{2,1}$ as well as a positive definite, diagonal matrix $\tSigma_{2,1}\in \mathbb R^{n_2, n_2}$.
\STATE Set $\mP_2 := \diag(\mV_{2,1}^H,\ \mU_{2,1}^H)$,\ $\mP:= \mP_2 \mP$.
%
\STATE Set
\[
\def\arraystretch{1.4}
\tJ := \mP_2\ \tJ\ \mP_2^H
=: \left[ \begin{array}{c|cc}
 \tJ_{1,1} & -\tJ_{2,1}^H & 0 \\
 \hline
 \tJ_{2,1} & \tJ_{2,2} & -\tJ_{3,2}^H \\
 0     & \tJ_{3,2} & \tJ_{3,3}
\end{array}\right],\qquad
\def\arraystretch{1}
\tR
:= \mP_2 \tR \mP_2^H
=: \left[ \begin{array}{c|cc}
\mR_1 & 0 & 0 \\
\hline
0 & 0 & 0 \\
0 & 0 & 0
\end{array}\right].
\]
(The lines indicate the partitioning of the block matrices~$\tJ$ and~$\tR$ in the previous step.)
\newline ----------- \textit{Step 2} -----------
\STATE i := 3
\WHILE{$n_{i-1} > 0$ \OR $\tJ_{i,i-1} \neq 0$}
\STATE Perform an SVD of $\tJ_{i,i-1}$ such that
\[
\tJ_{i,i-1}
= \mU_{i,i-1} \begin{bmatrix} \tSigma_{i,i-1} & 0\\ 0 & 0 \end{bmatrix} \mV^H_{i,i-1} ,
\]
with unitary matrices~$\mU_{i,i-1}$ and~$\mV_{i,i-1}$ as well as a positive definite, diagonal matrix $\tSigma_{i,i-1}\in \R^{n_i, n_i}$.
\STATE Set $\mP_i := \diag(\mI_{n_1},\ldots,\mI_{n_{i-2}},\ \mV_{i,i-1}^H,\ \mU_{i,i-1}^H )$, $\mP:= \mP_i \mP$.
%
\STATE Set
\[
\tJ
:= \mP_i\ \tJ\ \mP_i^H
=: \begin{bmatrix}
 \tJ_{1,1} & -\tJ_{2,1}^H & 0           &\cdots      & 0 \\
 \tJ_{2,1} & \tJ_{2,2}    & -\tJ_{3,2}^H &            & \vdots \\
          & \ddots      & \ddots      & \ddots     & \\
 \vdots   &             & \tJ_{i,i-1}  & \tJ_{i,i}   &-\tJ_{i+1,i}^H \\
    0     & \cdots      &     0       & \tJ_{i+1,i} & \tJ_{i+1,i+1}
\end{bmatrix} ,
\quad
\text{where } \tJ_{i,i-1} = [ \tSigma_{i,i-1}\quad 0].
\]
\STATE $i:= i+1$
\ENDWHILE
\newline ----------- \textit{Step 3} -----------
\STATE $s:= i$
\FOR{$i=2,\ldots,s$}
\STATE Set $\mJ_{i,i-1} :=\tJ_{i,i-1}$.
\ENDFOR
\FOR{$i=1,\ldots,s$}
\STATE Set $\mJ_{i,i} :=\tJ_{i,i}$.
\ENDFOR
\ENSURE Unitary matrix~$\mP$ satisfying~\eqref{matrices:staircase:J-R}.
\end{algorithmic}
\end{breakablealgorithm}
It is clear that Algorithm~\ref{algorithm:staircase:J-R} terminates after a finite number of steps either with~$n_{i-1}=0$ or~$\mJ_{i,i-1}=0$.
Note that the diagonal structure of the blocks~$\Sigma_{i,i-1}$ is destroyed in the next step of the algorithm except in the last step.
\end{proof}
Lemma~\ref{lem:SF} implies that if $\mR\geq 0$ then $\mR_1>0$ and $n_1=\rank \mR$.

Next, we give (a constructive proof of) a block decomposition for a matrix $\mA=\mJ-\mR$ with $\mR=\mR^H\geq 0$ and $\mJ=-\mJ^H$ under unitary congruence. This lemma will be used to construct a staircase form for matrix triples~$(\mE,\mJ,\mR)$ of semi-dissipative Hamiltonian DAE systems~\eqref{DAE:EA}.
\begin{lemma}[Unitary full rank decomposition for $\mJ-\mR$]
\label{lem:blocksvd}
Let $\mJ,\mR\in\Cnn$ satisfy $\mR=\mR^H\geq 0$ and $\mJ=-\mJ^H$,
then there exists a unitary matrix $\mP\in \Cnn$ such that
\begin{equation}\label{J-R:UFRD}
\mP\ (\mJ-\mR)\ \mP^H
=\begin{array}{l}
\begin{bmatrix} \mA_{1,1} & 0 \\ 0 & 0 \end{bmatrix}
\begin{array}{c} n_1 \\ n_2 \end{array} \\
  \quad n_1 \hspace{7pt} n_2
\end{array}
\end{equation}
with nonsingular $\mA_{1,1}\in\C^{n_1\times n_1}$ and $n_1,n_2\in\N_0$.
\end{lemma}
\begin{proof}[{Proof of Lemma~\ref{lem:blocksvd}}]
The proof is similar to that of Lemma~\ref{lem:SF} and given by the following constructive Algorithm~\ref{algorithm:UFRD:J-R}.
\begin{breakablealgorithm}
\caption{Unitary full rank decomposition for $\mJ-\mR$}
\label{algorithm:UFRD:J-R}
\begin{algorithmic}[1]
\REQUIRE $(\mJ,\mR)$
\newline ----------- \textit{Step 0} -----------
\STATE Perform a (spectral) decomposition of $\mR$ such that
\[
\mR
= \mU_R \begin{bmatrix} \tR_{1,1} & 0 \\ 0 & 0 \end{bmatrix} \mU_R^H ,
\]
with $\mU_R\in \Cnn$ unitary, $\tR_{1,1}=\tR_{1,1}^H\in\C^{\tilde n_1\times\tilde n_1}$ positive definite, and $\tilde n_1:=\rank\mR$.
\STATE Set $\mP := \mU^H_R$,
\[
\tJ
:= \mU^H_R\ \mJ\ \mU_R
=:\begin{bmatrix} \tJ_{1,1} & \tJ_{1,2}\\ \tJ_{2,1} & \tJ_{2,2} \end{bmatrix},\qquad
\tR:=\mU^H_R\ \mR\ \mU_R.
\]
\newline ----------- \textit{Step 1} -----------
\STATE Perform a (spectral) decomposition of $\tJ_{2,2}\in\C^{(n-\tilde n_1)\times(n-\tilde n_1)}$ such that
\[
\tJ_{2,2}
= \mU_{2} \begin{bmatrix} \tSigma_{2,2} & 0\\ 0 & 0 \end{bmatrix} \mU^H_{2},
\]
for some unitary matrix~$\mU_2$ and an invertible (square) matrix~$\tSigma_{2,2}\in\C^{\tilde n_2\times\tilde n_2}$.
\STATE Set $\mP_2 :=\diag(\mI, \mU_2^H)$,\ $\mP:= \mP_2 \mP$.
\STATE Set
\begin{equation} \label{lem:blocksvd:Step1}
\def\arraystretch{1.4}
\tJ := \mP_2\ \tJ\ \mP_2^H
=:
\left[ \begin{array}{c|cc}
 \tJ_{1,1} & -\tJ_{2,1}^H & -\tJ_{3,1}^H \\
 \hline
 \tJ_{2,1} & \tJ_{2,2} & 0 \\
 \tJ_{3,1} & 0 & 0
\end{array}\right]
\begin{array}{l}
 \tilde n_1 \\ \tilde n_2 \\ n-\tilde n_1 -\tilde n_2
\end{array},\
\def\arraystretch{1}
\tR
:= \mP_2\ \tR\ \mP_2^H
=:
\left[ \begin{array}{c|cc}
 \tR_{1,1} & 0 & 0\\
 \hline
 0 & 0 & 0\\
 0 & 0 & 0
\end{array}\right].
\end{equation}
(The lines indicate the partitioning of the block matrices~$\tJ$ and~$\tR$ in  Step~0.)
\newline ----------- \textit{Step 2} -----------
\STATE Perform an SVD of~$\tJ_{3,1}$ such that
\[
\tJ_{3,1}
=\mU_{3,1} \begin{bmatrix} \tSigma_{3,1} & 0\\ 0 & 0 \end{bmatrix} \mV^H_{3,1}\in\C^{(n-\tilde n_1 -\tilde n_2)\times\tilde n_1} ,
\]
with unitary matrices~$\mU_{3,1}$ and~$\mV_{3,1}$ as well as a nonsingular matrix~$\tSigma_{3,1}\in\C^{\tilde n_3\times \tilde n_3}$.
\STATE Set $\mP_3 := \diag(\mV_{3,1}^H,\ \mI,\ \mU_{3,1}^H)$,\ $\mP:= \mP_3 \mP$.
%
\STATE Set
\[
\def\arraystretch{1.4}
\tJ-\tR
:= \mP_3 (\tJ-\tR) \mP_3^H
=:
\left[ \begin{array}{cc|c|cc}
 \tJ_{1,1} -\tR_{1,1}&  -\tJ_{2,1}^H -\tR_{2,1}^H  & -\tJ_{3,1}^H & -\tJ_{4,1}^H & 0 \\
 \tJ_{2,1}-\tR_{2,1} & \tJ_{2,2} -\tR_{2,2}&  -\tJ_{3,2}^H  &  0 & 0 \\
 \hline
 \tJ_{3,1} & \tJ_{3,2} & \tJ_{3,3} & 0 & 0 \\
 \hline
 \tJ_{4,1} & 0 & 0 & 0 & 0 \\
 0 & 0 & 0 & 0 & 0
\end{array}\right]
\begin{array}{l}
 \tilde n_3 \\ \tilde n_1 -\tilde n_3 \\ \tilde n_2 \\ \tilde n_3 \\ \tilde n_4
\end{array}
\ ,
\]
with $\tJ_{4,1}=\tSigma_{3,1}$, $\tJ_{3,3}$ (in~\eqref{lem:blocksvd:Step1} it equals~$\tJ_{2,2}=\tSigma_{2,2}$) nonsingular and~$\tR_{2,2}$ (as principal submatrix of a positive definite matrix) positive definite.
(The lines indicate the partitioning of the block matrices~$\tJ$ and~$\tR$ in~\eqref{lem:blocksvd:Step1}.)
\ENSURE Unitary matrix~$\mP$ satisfying~\eqref{J-R:UFRD}.
\end{algorithmic}
\end{breakablealgorithm}

Performing block Gaussian elimination we see immediately that the block consisting of the first four rows and columns in~$\tJ-\tR$ is square and invertible and corresponds to the matrix~$\mA_{1,1}$ in the assertion, i.e. $n_1 =\tilde n_1 +\tilde n_2 +\tilde n_3$.
\end{proof}

\subsection{Hypocoercive matrices}\label{sec:hypoODE}
In this section we recall the concept of hypocoercivity for linear ODEs and relate it to the staircase form of the last subsection.
\begin{definition}[{\cite{AAC20}}] \label{def:matrix:hypocoercive}
A matrix $\mC\in\Cnn$ is called~\emph{coercive} if its Hermitian part $\mCH$ is positive definite, and it is called~\emph{hypocoercive} if the spectrum of~$\mC$ lies in the \emph{open} right half plane.

For practical reasons, a matrix~$\mA\in\Cnn$ is called \emph{negative hypocoercive} if the spectrum of~$\mA$ lies in the \emph{open} left half plane.
\end{definition}
Hypocoercive matrices are often called~\emph{positively stable}, whereas negative hypocoercive matrices are often called~\emph{stable}.
We use the notion of hypocoercivity to emphasize the analogous situation in partial differential equations, see~\cite{AAC18,ArEr14,Vi09}.

To decide if a matrix with positive semi-definite Hermitian part is hypocoercive, or equivalently, if a semi-dissipative matrix is negative hypocoercive, one can proceed as follows:
\begin{proposition}[Lemma 3.1 in~\cite{MMS16}, Lemma~2.4 in~\cite{AAC18}]
\label{prop:border}
Let $\mJ,\mR\in\Cnn$ be such that $\mJ^H =-\mJ$ and $\mR^H =\mR\geq 0$.
Then, $\mJ+\mR$ (resp.~$\mJ-\mR$) has an eigenvalue on the imaginary axis if and only if $\mR v =0$ for some eigenvector~$v$ of~$\mJ$. 
\end{proposition}
Note that, due to the assumptions, purely imaginary eigenvalues of $\mJ+\mR$ are necessarily semi-simple, see also~\cite{MehMW18,MehMW20}.
Therefore, a matrix~$\mC$ with positive semi-definite Hermitian part is hypocoercive if and only if no eigenvector of the skew-Hermitian part~$\mCS$ lies in the kernel of the Hermitian part~$\mCH$.
The latter condition is well known in control theory, and equivalent to the following statements:
\begin{lemma} \label{lem:Equivalence}
Let $\mJ,\mR\in\Cnn$ satisfying $\mR=\mR^H\geq 0$ and $\mJ=-\mJ^H$. 
Then the following conditions are equivalent:
\begin{enumerate}[(B1)]
\item \label{B:KRC}
There exists $m\in\N_0$ such that
\begin{equation}\label{condition:KalmanRank}
 \rank[{\mR},\mJ{\mR},\ldots,\mJ^m {\mR}]=n \,.
\end{equation}
\item \label{B:Tm}
There exists $m\in\N_0$ such that
\begin{equation}\label{Tm:J-R}
 \mT_m :=\sum_{j=0}^m \mJ^j \mR (\mJ^H)^j > 0 \,.
\end{equation}
%
\item \label{B:PBH:EVec}
No eigenvector of~$\mJ$ lies in the kernel of~$\mR$.
\item \label{B:PBH:EVal}
$\rank [\lambda \mI-\mJ, \mR] =n$ for every $\lambda \in \mathbb C$ , in particular for every eigenvalue~$\lambda$ of~$\mJ$.
\end{enumerate}
Moreover, the smallest possible~$m\in\N_0$ in~\ref{B:KRC} 
and~\ref{B:Tm} coincide. 
\end{lemma}
\begin{proof}
The statement of Lemma~\ref{lem:Equivalence} and its proof are classical, see e.g. Theorem~6.2.1 of~\cite{Da04} for real matrices and Proposition~1 of~\cite{AAC18} for complex matrices.
\end{proof}
To summarize, a matrix with positive semi-definite Hermitian part is hypocoercive if one (hence, all) conditions in Lemma~\ref{lem:Equivalence} are satisfied.
\new{The relationship between conditions~\ref{B:PBH:EVec}--\ref{B:PBH:EVal} and stability concepts for semigroups is also known in the infinite-dimensional case, see Theorem~14 in~\cite{BaVu90}.
Moreover, an alternative formulation of the Kalman rank condition~\ref{B:KRC} which is applicable in the infinite-dimensional setting (and allows to characterize the hypocoercivity index) is given in Remark~4 of~\cite{AAC18}.
In a forthcoming work~\cite{AAM21} we develop our concepts in the infinite-dimensional setting including a possible extension of Lemma~\ref{lem:Equivalence}.
}

Next, we define the hypocoercivity index for matrices with positive semi-definite Hermitian part as in~\ref{B:Tm}:
\begin{definition}[{\cite{AAC20}}]
Let $\mJ,\mR\in\Cnn$ satisfy $\mR=\mR^H\geq 0$ and $\mJ=-\mJ^H$.
The~\emph{hypocoercivity index (HC-index)~$m_{HC}$ of the matrix~$\mJ+\mR$} is defined as the smallest integer~$m\in\N_0$ such that~\eqref{Tm:J-R} holds.
For matrices $\mJ+\mR$ that are not hypocoercive we set $m_{HC}=\infty$.

For practical reasons, we define the HC-index~$\mHC$ also for semi-dissipative matrices~$\mJ-\mR$ as the smallest integer~$m\in\N_0$ such that~\eqref{Tm:J-R} holds.
\end{definition}
Note that a hypocoercive matrix~$\mC=\mJ+\mR$ is coercive if and only if $m_{HC}=0$.
Similarly, a semi-dissipative matrix~$\mA=\mJ-\mR$ is dissipative if and only if $m_{HC}=0$.

It is obvious that the HC-indices of $\mJ\pm\mR$ and of a unitarily congruent matrix~$\mP(\mJ\pm\mR)\mP^H$ coincide.

\begin{lemma} \label{lem:SF:hypocoercive}
Let~$\mA =\mJ -\mR$ be a semi-dissipative matrix with skew-Hermitian matrix~$\mJ =\mAS$ and Hermitian matrix~$\mR =-\mAH$.
Then there exists a unitary matrix~$\mP\in\Cnn$ such that~\eqref{matrices:staircase:J-R} holds, where $\mR_1>0$ is a positive definite Hermitian matrix.
Moreover, the matrix $\mA$ is negative hypocoercive if and only if~$n_s =0$, i.e., the last row and last column in~$\mP\mJ\mP^H$ of~\eqref{matrices:staircase:J-R} are absent, and the HC-index of~$\mA$ is~$\mHC =s-2$.
\end{lemma}
\begin{proof}
The first statement follows from Lemma~\ref{lem:SF}.

Due to Lemma~\ref{lem:Equivalence}, the (hypocoercivity) condition~\eqref{Tm:J-R} is equivalent to~\ref{B:PBH:EVal}. 
Using again the staircase form~\eqref{matrices:staircase:J-R} we observe that~\ref{B:PBH:EVal} holds if and only if $n_s=0$.
In this case the value of the HC-index $\mHC=s-2$ can be deduced from the staircase form~\eqref{matrices:staircase:J-R} and the rank condition~\ref{B:KRC}.
This finishes the proof (of the second statement).
\end{proof}
\begin{remark}\label{rem:controlTheory}
The characterization of (negative) hypocoercive matrices is also related to results in control theory:
Consider a state-space system
\begin{equation}\label{ODE:AxBu}
 \dot x(t) =\mA x +\mB u
\end{equation}
for constant matrices $\mA,\mB\in\Cnn$. 
\begin{itemize}
\item
A pair~$(\mA,\mB)$ of square matrices~$\mA,\mB\in\Cnn$ is called \emph{controllable} if the controllability matrix~$[\mB, \mA\mB, \mA^2 \mB,\ldots,\mA^{n-1} \mB]$ has full rank.
For a controllable pair~$(\mA,\mB)$, the smallest possible integer~$k$ such that the controllability (sub)matrix $[\mB, \mA\mB, \mA^2 \mB,\ldots,\mA^{k-1} \mB]$ has full rank, is called the \emph{controllability index}, see e.g.~page 121 of~\cite{Wo85} or Paragraph 10.1.4.6 of~\cite{Le11}.
Thus, for a semi-dissipative matrix~$\mJ-\mR$, its HC-index is one less than the controllability index of~$(\mA,\mB)=(\mJ,\mR)$.
Using this interpretation, Lemma~\ref{lem:SF:hypocoercive} is a special case of results in control theory, see e.g.~Theorem 6.7.1 of~\cite{Da04}.
\item
A state-space system~\eqref{ODE:AxBu} is called ~\emph{asymptotically controllable to~$0$} if for all $x_0\in\Cn$ there exists a control~$u(t)$ such that the solution~$x(t)$ of~\eqref{ODE:AxBu} with $x(0)=x_0$ satisfies~$\lim_{t\to\infty} x(t) =0$, see Section~5.5 in~\cite{So98}.
For a semi-dissipative matrix~$\mJ-\mR$, the state-space system~\eqref{ODE:AxBu} with~$(\mA,\mB)=(-\mR,\mJ)$ is asymptotically controllable to~0 if~$\mJ-\mR$ has a finite HC-index (e.g. take $u=x$).
However, the converse is not true: Consider the semi-dissipative matrix $\mJ-\mR=\diag(i,-1)$ with HC-index $\mHC=\infty$ which is asymptotically controllable to~0 using $u=i\ x$.
\end{itemize}
\end{remark}

With an eye towards numerical computations we shall finally discuss perturbations of (negative) hypocoercive matrices.
An interesting open problem is the question of the smallest perturbation that either increases the HC-index or even makes a (negative) hypocoercive matrix not (negative) hypocoercive, which by the previous analysis is equal to the distance to uncontrollability for structured pairs $(\mJ,\mR)\in i\Hn\times\PSDHn$.

To study this question, it is reasonable to preserve the structure of matrices (or systems) in the perturbation analysis, such as for semi-dissipative Hamiltonian ODE systems.
The set of semi-dissipative matrices is convex, i.e. for semi-dissipative matrices $\mA=\mJ-\mR$, $\mA_1=\mJ_1-\mR_1$ and non-negative constant $\delta\geq0$, the matrices $\mA+\delta\mA_1=(\mJ+\delta\mJ_1) -(\mR+\delta\mR_1)$ are again semi-dissipative (with Hermitian part $\mR+\delta\mR_1$ and skew-Hermitian part~$\mJ+\delta\mJ_1$).
Moreover, it is reasonable to restrict to perturbations which do not change the rank of the Hermitian part, i.e.
\begin{equation}\label{perturbation:admissible}
\rank\mR
=\rank(\mR+\delta\mR_1)
\quad\text{for all } \delta\geq 0.
\end{equation}
\begin{eexample}\label{ex:HCmatrix+perturbation}
Consider the semi-dissipative matrix~$\mA=\mJ-\mR$ and perturbation matrices $\mJ_1,\mR_1$,
\begin{align*}
\mJ
&:=
\begin{bmatrix}
 0 & 1 & 0 & 0 \\
 -1 & 0 & 1 & 0 \\
 0 & -1 & 0 & 0 \\
 0 & 0 & 0 & 0
\end{bmatrix} ,
& %
\mR
&:=
\begin{bmatrix}
 0 & 0 & 0 & 0 \\
 0 & 0 & 0 & 0 \\
 0 & 0 & 1 & 0 \\
 0 & 0 & 0 & 1
\end{bmatrix} ,
& %
\mJ_1
&:=
\begin{bmatrix}
 0 & 0 & 0 & 1 \\
 0 & 0 & 0 & 0 \\
 0 & 0 & 0 & 0 \\
 -1 & 0 & 0 & 0
\end{bmatrix} ,
& %
\mR_1
:=
\begin{bmatrix}
 0 & 0 & 0 & 0 \\
 0 & 0 & 0 & 0 \\
 0 & 0 & 0 & 0 \\
 0 & 0 & 0 & 0
\end{bmatrix} ,
\end{align*}
such that~\eqref{perturbation:admissible} holds.
For~$\delta\geq0$, the HC-index of the matrices $\mA+\delta\mA_1=(\mJ+\delta\mJ_1) -(\mR+\delta\mR_1)$ satisfies
\[
 \mHC = \begin{cases}
         1 & \text{if } \delta>0 , \\
         2 & \text{if } \delta=0 .
        \end{cases}
\]
\end{eexample}
Example~\ref{ex:HCmatrix+perturbation} shows that it is easy to decrease the HC-index with arbitrary small perturbations which preserve the structure.
By contrast, the hypocoercivity condition~\eqref{Tm:J-R} shows that small enough perturbations cannot increase the HC-index.

\subsection{Hypocoercivity and Lyapunov stability}\label{ssec:ODE:lyastab}

It is well known that an ODE~\eqref{ODE:A} is \emph{(Lyapunov) stable} if all eigenvalues of~$\mA$ have non-positive real part and the eigenvalues on the imaginary axis are semi-simple, and~\eqref{ODE:A} is \emph{asymptotically stable} if all eigenvalues of~$\mA$ have negative real part.

A semi-dissipative Hamiltonian ODE~\eqref{ODE:A} with $\mA=\mJ-\mR$ is (Lyapunov) stable, since for all solutions~$x(t)$ of~\eqref{ODE:A} the Euclidean norm is nonincreasing
\begin{equation} \label{energy:estimate}
\ddt \|x(t)\|^2
= \ip{\mA x(t)}{x(t)} +\ip{x(t)}{\mA x(t)}
= \ip{x(t)}{(\mA^H +\mA)x(t)}
\leq 0 ,
\quad t\geq 0 .
\end{equation}
Due to Proposition~\ref{prop:border} and Lemma~\ref{lem:Equivalence},
a semi-dissipative Hamiltonian ODE~\eqref{ODE:A} is asymptotically stable if and only if the HC-index (of its system matrix~$\mA$) is finite.

Phenomenologically, the HC-index describes the structural complexity of the intertwining of the two matrices $\mR$ and $\mJ$ (see~\cite{AAC18} for illustrating examples).
Moreover, for a semi-dissipative Hamiltonian ODE~\eqref{ODE:A}, the HC-index characterizes the decay of its \emph{propagator norm} for short time.
We denote the solution semigroup pertaining to~\eqref{ODE:A} by $S(t):=e^{\mA t}\in \Cnn$.
The \emph{short-time decay} of its spectral norm~$\|S(t)\|_2$ is related to the HC-index as follows:
\begin{theorem}[\cite{AAC20}]\label{th:HC-decay}
Consider a semi-dissipative Hamiltonian ODE~\eqref{ODE:A} whose system matrix~$\mA$ has finite HC-index.
Its (finite) HC-index is $\mHC\in\N_0$ if and only if
\begin{equation}\label{short-t-decay}
\|S(t)\|_2
=
1 -ct^a +\mathcal O(t^{a+1})
\quad\text{for } t\to0+\,,
\end{equation}
where $c>0$ and $a=2\mHC+1$.
\end{theorem}
\begin{eexample}[Example~5.2 in~\cite{AAS19}] \label{ex:envelope}
We consider ODE~\eqref{ODE:A} with the semi-dissipative matrix
\begin{equation}\label{ODE:A:envelope}
 \mA :=\begin{bmatrix} -1 & 1 \\ -1 & 0 \end{bmatrix} ,
\end{equation}
whose eigenvalues are $\lambda_\pm =(-1\pm i\sqrt{3})/2$.
The Hermitian part $\mAH =\diag({-1}, 0)$ is only semi-definite
and the semi-dissipative matrix~$\mA$ has HC-index $\mHC=1$.
The squared \emph{propagator norm}, see~\cite{ASS20,GaMi13}, satisfies
\begin{equation} \label{ex:norm2:S}
\begin{split}
 \|e^{\mA t}\|^2_2
 &= \tfrac16 \ \bigg(\Big(\sqrt{-2\cos(\sqrt{3}t) + 14} + \sqrt{-2\cos(\sqrt{3}t) + 2}\Big)\ \sqrt{-2\cos(\sqrt{3}t) + 2} + 6 \bigg) e^{-t} \\
 &\sim 1 -\tfrac16 t^3 +\bigO(t^4) \quad\text{for } t\to 0+ ,
\end{split}
\end{equation}
which illustrates the result of Theorem~\ref{th:HC-decay} with~$a =2\mHC+1 =3$.
The kernel of~$\mAH$ is one-dimensional and it is spanned by the normalized vector $\bx_0 =[0, 1]^\top$.
The squared norm of the solution of~\eqref{ODE:A} with initial condition $\bx(0)=\bx_0$ is given by
\begin{equation} \label{ex:norm2:x}
\begin{split}
\|\bx(t)\|^2_2
&=
\tfrac19  e^{-t}\ {\bigg(\sqrt{3} \sin\big(\tfrac{\sqrt{3}}{2} \ t\big) - 3 \, \cos\big(\tfrac{\sqrt{3}}{2} \ t\big)\bigg)}^{2} + \tfrac43 \, e^{-t} \sin\big(\tfrac{\sqrt{3}}{2} \ t\big)^{2} \\
&\sim 
1 -\tfrac23 t^3 +\bigO(t^4) 
\quad\text{for } t\to 0+.
\end{split}
\end{equation}
Thus the Taylor expansions of~\eqref{ex:norm2:S} and~\eqref{ex:norm2:x} have the same form $1-c t^3 +\bigO(t^4)$ as $t\to0+$, whereas the propagator norm decays slightly slower than the solution starting at the vector $\bx_0$, 
see also Figure~\ref{fig:envelope}.
\begin{figure}[h!]
 \centering
 \includegraphics[width=0.75\textwidth]{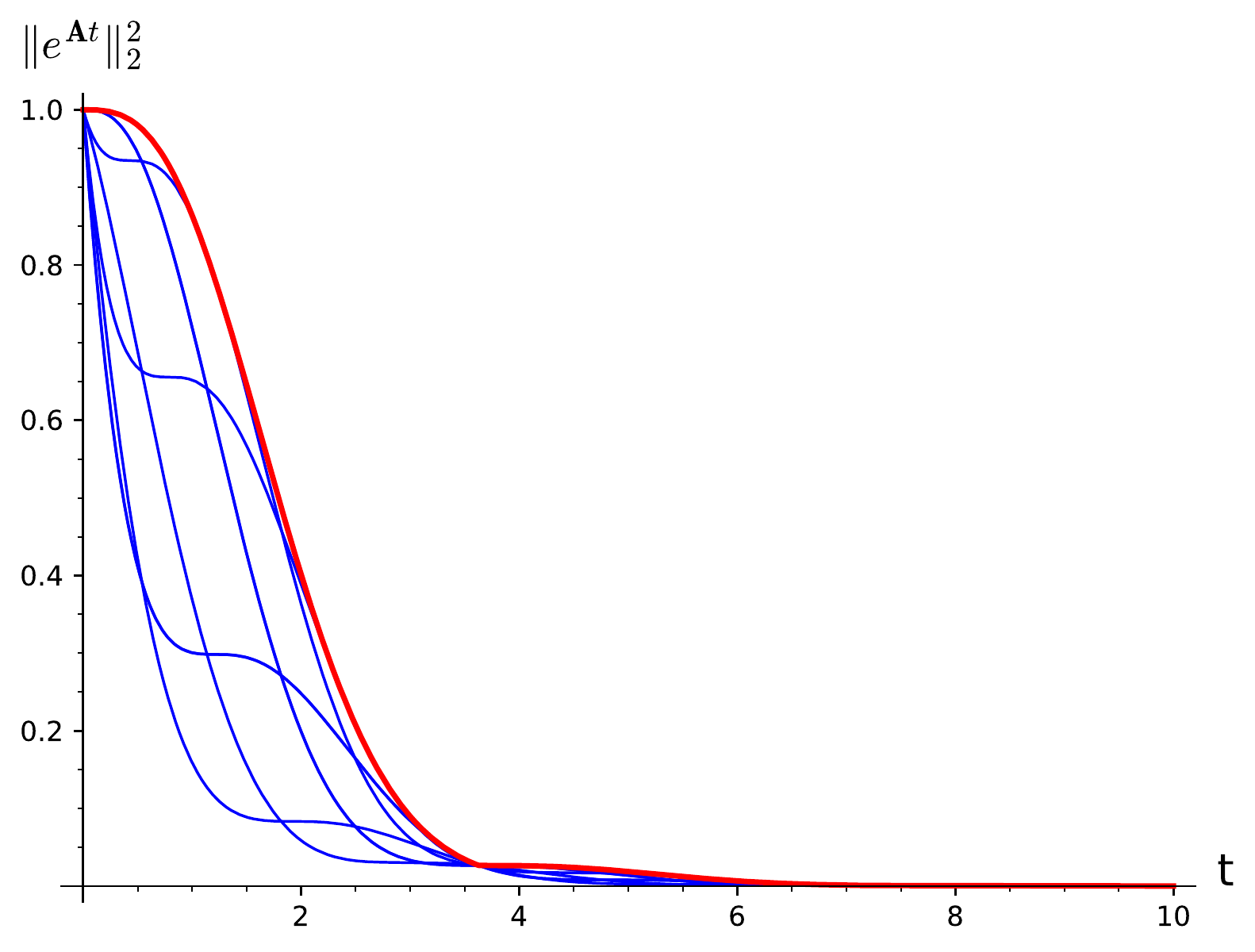}
 \caption{For the ODE~\eqref{ODE:A} with $\mA=\left[ \begin{array}{cc} -1 & 1 \\ -1 & 0 \end{array}\right]$, the squared propagator norm (red line) and the squared norms of a family of solutions of~\eqref{ODE:A} (blue lines) are plotted. The squared propagator norm satisfies~$\|e^{\mA t}\|^2_2 \sim 1 -t^3/6 +\bigO(t^4)$ for $t\to 0+$, and it is not continuously differentiable at $t=2\pi/\sqrt3$.
 Moreover, it is the envelope of $\|x(t)\|_2^2$ for all solutions with~$\|x(0)\|_2^2=1$.}
 \label{fig:envelope}
\end{figure}
\end{eexample}

\medskip
The stability of~\eqref{ODE:A} can also be characterized via (solutions of) Lyapunov matrix inequalities, see e.g.~\cite{St75}:
\begin{proposition}[Lyapunov] \label{prop:stability:LF} 
Consider a system of ODEs~\eqref{ODE:A} with~$\mA\in\Cnn$.
\begin{enumerate}[(S1)]
\item \label{P:stable}
System~\eqref{ODE:A} is stable 
if and only if there exists~$\mX\in\PDHn$ such that~$\mA^H \mX+\mX\mA \leq 0$.
\item \label{P:asymptotically-stable}
System~\eqref{ODE:A} is asymptotically stable 
if and only if there exists~$\mX\in\PDHn$ such that~$\mA^H \mX+\mX\mA < 0$.
\item
System~\eqref{ODE:A} is unstable if there exists a Hermitian matrix~$\mX\in\Cnn$ that is either negative definite or indefinite such that~$\mA^H \mX+\mX\mA <0$.
\end{enumerate}
\end{proposition}
%
The characterization of stability in Proposition~\ref{prop:stability:LF} allows to define quadratic Lyapunov functionals as a squared weighted norm $\|\cdot\|_\mX^2 =\ip{\cdot}{\mX\cdot}$ for some appropriate matrix~$\mX\in\PDHn$.
For semi-dissipative Hamiltonian ODEs~\eqref{ODE:A} with matrix~$\mA$, the identity matrix~$\mX =\mI$ satisfies~$\mA^H \mX+\mX\mA \leq 0$, hence the Euclidean norm is a (non-strict) Lyapunov functional, see~\eqref{energy:estimate}.
If~\eqref{ODE:A} is even asymptotically stable, then the strict inequality of~\ref{P:asymptotically-stable} can be improved to the form of~\eqref{example:X:LMI}, see below.
This then yields exponential decay of the (strict) Lyapunov functional $\|\cdot\|_\mX^2$ along solutions of~\eqref{ODE:A}.

A negative hypocoercive (or stable) matrix~$\mA$ does not necessarily have a negative definite Hermitian part:
\begin{eexample}[Example~5.3 in~\cite{AAS19}] \label{example:AH:indefinite}
The matrix
\begin{equation} \label{example:A}
\mA
:= \begin{bmatrix} -19 & 6 \\ -6 & 1 \end{bmatrix}
\end{equation}
is negative hypocoercive (or stable) since its eigenvalues are $\lambda_1=-1$ and $\lambda_2=-17$.
However, its Hermitian part~$\mAH=\diag(-19,1)$ is indefinite.
Nonetheless, it is possible to construct matrices~$\mX\in\PDHn$ which satisfy the Lyapunov matrix inequality $\mA^H \mX+\mX\mA <0$.
In fact, the family of Hermitian matrices
\begin{equation} \label{example:X}
\mX
:=
\begin{bmatrix} 1 & 3 \\ -3 & -1 \end{bmatrix}
\begin{bmatrix} b_1 & 0 \\ 0 & b_2 \end{bmatrix}
\begin{bmatrix} 1 & -3 \\ 3 & -1 \end{bmatrix}
\qquad \text{with } b_1,b_2>0 \ ,
\end{equation}
satisfies
\begin{equation} \label{example:X:LMI}
\mA^H \mX+\mX\mA
\leq 2\lambda_1\mX <0 \ .
\end{equation}
\end{eexample}

\bigskip
Another point of view is that for a stable system of ODEs~\eqref{ODE:A} there exists a representation as a semi-dissipative Hamiltonian DAE~\eqref{DAE:EA} with positive definite Hermitian matrix~$\mE$:
More precisely, an ODE~\eqref{ODE:A} with matrix~$\mA$ is stable (resp.~asymptotically stable) if and only if there exists a semi-dissipative (resp.~dissipative) Hamiltonian DAE~\eqref{DAE:EA} with pencil
\begin{equation}\label{matrix:A:dH}
\lambda \mE -(\mJ-\mR)
\quad\text{such that }
\mE =\mX \text{ and } \mX\mA =\mJ-\mR
\end{equation}
for some~$\mX,\mJ,\mR\in\Cnn$ satisfying $\mX=\mX^H >0$, $\mR=\mR^H\geq 0$ (resp. $>0$) and $\mJ=-\mJ^H$.
This result is well known in a much more general context, see e.g.~\cite{Wil72a,Wil72b}, and~\cite{BeMeVD19} in the context of port-Hamiltonian systems.

In the case of semi-dissipative Hamiltonian DAEs~\eqref{DAE:EA} with singular~$\mE$, the question of stability is much more complex (than for DAEs~\eqref{DAE:EA} with regular~$\mE$), since the classical relation between the existence of positive definite solutions to Lyapunov equations and stability of a system does not hold any longer, see~\cite{Sty02,Sty02a} for a detailed analysis.
The essential difference is that the solution of the Lyapunov equation need not be semi-definite, only the part associated with the dynamic part, and the right hand side has to be adapted.

\subsection{Construction of strict Lyapunov functionals}
\label{sec:LyaFunc}

\new{Consider a semi-dissipative Hamiltonian ODE~\eqref{ODE:A} with matrix~$\mA$. 
Using the full spectral information about the eigenvalue/eigenvector structure of~$\mA$ allows to construct strict Lyapunov functionals which provide the sharp decay rate $\lambda$ in~\eqref{exp-decay} (for large time), see~\cite{AAC16,AASt15,ArEr14,ArJiWo19}.
The aim of Algorithm~\eqref{algorithm:1Pi} below is to construct a strict Lyapunov functional without the full eigenvalue/eigenvector information.
Therefore, it is expected to yield only a suboptimal decay rate.
However, the construction can be applied also in the infinite-dimensional setting, see~\cite{AAC18,AAC20,AAM21}.} 

\smallskip
The HC-index, resp. controllability index appears in the construction of strict Lyapunov functionals for semi-dissipative \new{Hamiltonian} ODEs, see Section~4 of~\cite{AAC20}, resp. Definition~6.3.2 of~\cite{Mi09thesis}.
In~\S2.3 of~\cite{AAC18}, it was used to give efficient guidelines for this construction in cases with $\dim(\ker\mR)$ small.
A recent improvement allows to construct a strict Lyapunov functional for all semi-dissipative \new{Hamiltonian} ODEs with finite HC-index~$\mHC$, namely in~$\mHC$ steps, see~\cite{AAC20}.

\new{First, we recall the construction of an ansatz for strict Lyapunov functionals in Algorithm~\ref{algorithm:1Pi}, which is taken from~\cite{AAC20}.
There we derive explicit restrictions on~$\eps_j$ (relative to other parameters) such that a suitable choice of $\eps_j$ turns the ansatz in step~10 of Algorithm~\ref{algorithm:1Pi} into a strict Lyapunov functional.
Then, we will show that Algorithm~\ref{algorithm:1Pi} uses implicitly a staircase form (similar to the one in Lemma~\ref{lem:SF}).
} 

Consider a semi-dissipative matrix~$\mA=\mJ-\mR$ with finite HC-index, then our algorithm reads as follows:

\begin{breakablealgorithm}
\caption{Construction of a strict Lyapunov functional}
\label{algorithm:1Pi}
\begin{algorithmic}[1]
\REQUIRE $\Pi_0:=\mI$, $\tA_0:=-\mJ$, $\tB_0:=\mR$, $j:=1$
\STATE Construct an orthogonal projection~$\widetilde{\Pi}_j$ onto $\ker \big(\tB_{j-1} \tB_{j-1}^H\big)$.
\STATE $\Pi_j :=\widetilde{\Pi}_j \Pi_{j-1}$
\WHILE{$\Pi_j \ne0$}
\STATE Set $\tA_j :=\Pi_j \tA_{j-1} \Pi_j$, $\tB_j :=\Pi_j \tA_{j-1} (\Pi_{j-1} -\Pi_j)$.
\STATE $j := j+1$
\STATE Construct an orthogonal projection~$\widetilde{\Pi}_j$ onto $\ker \big(\tB_{j-1} \tB_{j-1}^H\big)$.
\STATE $\Pi_j :=\widetilde{\Pi}_j \Pi_{j-1}$
\ENDWHILE
\STATE $\mHC :=j-1$
\STATE Set $\mX :=\Pi_0 +\sum_{j=1}^{\mHC} \eps_j \big(\tA_{j-1}\Pi_j +\Pi_j \tA_{j-1}^H\big)$ for sufficiently small $\eps_j>0$.
\ENSURE $\|\cdot\|_\mX^2 :=\ip{\cdot}{\mX\cdot}$
\end{algorithmic}
\end{breakablealgorithm}
\begin{remark}\label{rem:large}\new{
Algorithm~\ref{algorithm:1Pi} can be easily implemented as a numerically stable algorithm for small and medium size problems via $QR$ and singular value decompositions. 
The most challenging part is the number of consecutive rank decisions, which is critical in this context. 
However, this has been dealt with in staircase algorithms and actually is easier in the case of semi-dissipative \new{Hamiltonian} ODEs and also DAEs, \cite{BruM07,ByeMX07}. 
For large scale problems it is possible to use Krylov subspace methods  and truncated singular value or skinny $QR$ decompositions, see \cite{scalapack}. 
In view of the applications from discretized partial differential equations with constraints, these kernels are often available directly from the modeling or finite element discretization, see e.g. \cite{EmmM13,GraMQSW16}.
}
\end{remark}
To illustrate the construction of \new{(an ansatz for)} a strict Lyapunov functional as a squared weighted norm~$\|\cdot\|_\mX^2$ for a semi-dissipative matrix~$\mA$ with finite HC-index, we consider the following example inspired by the analysis of a kinetic transport-reaction model in one space dimension~\cite{AAC18}:
\begin{eexample}
Consider the matrix~$\mA=\mJ-\mR$ with
\begin{align} \label{ex:BGK:JR}
\mJ
&:= -i
\begin{bmatrix}
0 & \sqrt 1 & 0 & 0 & 0 \\
\sqrt 1 & 0 & \sqrt 2 & 0 & 0 \\
0 & \sqrt 2 & 0 & \sqrt 3 & 0 \\
0 & 0 & \sqrt 3 & 0 & \sqrt 4 \\
0 & 0 & 0 & \sqrt 4 & 0
\end{bmatrix} \,,
& %
\mR
&:= \diag(0,0,0,1,1) \,.
\intertext{The semi-dissipative matrix $\mA$ has HC-index~$3$.
\newline\indent
Define $\Pi_0:=\mI$, $\tA_0:=-\mJ$, $\tB_0:=\mR$, $j:=1$.
The orthogonal projection $\widetilde{\Pi}_1$ onto $\ker \big(\tB_0 \tB_0^H\big)$ is given by $\widetilde{\Pi}_1 =\diag(1,1,1,0,0)$.
Define $\Pi_1:=\widetilde{\Pi}_1 \Pi_0 =\widetilde{\Pi}_1$.
Then, the matrix $\tA_0$ can be written as $\tA_0 =\tA_1 +\tB_1 +\tX_1 +\tY_1$ with
}
\tA_1
&:= \Pi_1 \tA_0 \Pi_1\,,
&
\tB_1
&:= \Pi_1 \tA_0 (\mI -\Pi_1) \,,
\nonumber \\
\tX_1
&:= (\mI -\Pi_1) \tA_0 \Pi_1 \,,
&
\tY_1
&:= (\mI -\Pi_1) \tA_0 (\mI -\Pi_1) \,.
\nonumber
\intertext{In our example, we compute}
\tA_1
&=\Pi_1 \tA_0 \Pi_1
= i\ \left[ \begin{array}{ccc|cc}
0 & \sqrt 1 & 0 & 0 & 0 \\
\sqrt 1 & 0 & \sqrt 2 & 0 & 0 \\
0 & \sqrt 2 & 0 & 0 & 0 \\
\hline
0 & 0 & 0 & 0 & 0 \\
0 & 0 & 0 & 0 & 0
\end{array}\right] ,
& %
\tB_1
&= i\ \left[\begin{array}{ccc|cc}
0 & 0 & 0 & 0 & 0 \\
0 & 0 & 0 & 0 & 0 \\
0 & 0 & 0 & \sqrt 3 & 0 \\
\hline
0 & 0 & 0 & 0 & 0 \\
0 & 0 & 0 & 0 & 0
\end{array}\right] .
\nonumber
\intertext{
The projection $\widetilde{\Pi}_2$ onto $\ker \big(\tB_1 \tB_1^H\big) =\ker(\diag(0,0,3,0,0))$ is given by $\widetilde{\Pi}_2 =\diag(1,1,0,1,1)$.
Therefore, $\Pi_2 =\widetilde{\Pi}_2 \Pi_1 =\diag(1,1,0,0,0)$.
Since $\Pi_2\ne 0$, we continue the iteration on the upper left $3\times 3$ block of~$\tA_1$ (indicated by the partition lines above).
We define
}
\def\arraystretch{1.2}
\tA_2
&:= \Pi_2 \tA_1 \Pi_2
= i\ \left[\begin{array}{ccc|cc}
0 & \sqrt 1 & 0 & 0 & 0 \\
\sqrt 1 & 0 & 0 & 0 & 0 \\
0 & 0 & 0 & 0 & 0 \\
\hline
0 & 0 & 0 & 0 & 0 \\
0 & 0 & 0 & 0 & 0
\end{array}\right] ,
& %
\tB_2
&:= i\ \left[\begin{array}{ccc|cc}
0 & 0 & 0 & 0 & 0 \\
0 & 0 & \sqrt 2 & 0 & 0 \\
0 & 0 & 0 & 0 & 0 \\
\hline
0 & 0 & 0 & 0 & 0 \\
0 & 0 & 0 & 0 & 0
\end{array}\right] ,
\nonumber
\intertext{
such that $\tB_2 \tB_2^H =\diag(0,2,0,0,0)$, $\widetilde{\Pi}_3 =\diag(1,0,1,1,1)$, and $\Pi_3 :=\widetilde{\Pi}_3 \Pi_2 =\diag(1,0,0,0,0)$.
Since $\Pi_3\ne 0$, we continue the iteration once more and define
}
\tA_3
&:= \Pi_3 \tA_2 \Pi_3 =0 \,, \qquad
& %
\tB_3
&:= i\ \begin{bmatrix}
0 & 1 & 0 & 0 & 0 \\
0 & 0 & 0 & 0 & 0 \\
0 & 0 & 0 & 0 & 0 \\
0 & 0 & 0 & 0 & 0 \\
0 & 0 & 0 & 0 & 0
\end{bmatrix} \,,
\nonumber
\end{align}
such that $\tB_3 \tB_3^H =\diag(1,0,0,0,0)$ and $\widetilde{\Pi}_4 =\diag(0,1,1,1,1)$.
Therefore, $\Pi_4 =\widetilde{\Pi}_4 \Pi_3 =0$ 
which indicates the end of the iteration.
Finally, for sufficiently small~$\eps_j>0$, the ansatz for a (strict) Lyapunov functional is given by $\|\cdot\|_\mX^2$ with
\begin{equation} \label{matrix:P}
\mX
:=\mI +\sum_{j=1}^{\mHC} \eps_j (\tA_{j-1}\Pi_j +\Pi_j \tA_{j-1}^H)
=\begin{bmatrix}
1 & -i \eps_3 & 0 & 0 & 0 \\
i \eps_3 & 1 & -i \eps_2 \sqrt 2 & 0 & 0 \\
0 & i \eps_2 \sqrt 2 & 1 & -i \eps_1 \sqrt 3 & 0 \\
0 & 0 & i \eps_1 \sqrt 3 & 1 & 0 \\
0 & 0 & 0 & 0 & 1
\end{bmatrix} \,,
\end{equation}
where $\mHC=3$.
In this way we recover the ansatz for a (strict) Lyapunov functional following our guidelines for ODE systems exhibiting hypocoercivity in~\S2.2.3 of~\cite{AAC18}, see also the ansatz for a linearized BGK model on a one-dimensional torus in Equation (69) of~\cite{AAC18}.
\end{eexample}

\subsubsection*{An alternative staircase form}

The matrix pair~$(\mJ,\mR)$ in~\eqref{ex:BGK:JR} is already in a staircase form where the order of the basis is reversed in comparison with~\eqref{matrices:staircase:J-R}.
Such an alternative staircase form has been considered in~\cite{vD81}, where for unstructured matrix pairs~$(\mA,\mB)$, a unitary matrix~$\mT$ is constructed such that the transformed matrix pencil~$[\lambda \mI -\mA || \mB]$ \new{has a block matrix structure of the form}
\begin{equation} \label{staircase:VD}
\begin{split}
& 
[ \mT^{-1}(\lambda \mI -\mA)\mT || \mT^{-1} \mB]
\\
&=:
\left[\begin{array}{c|cccccc||c}
\lambda \mI_{\tau_k} -\mA_k & 0 & & \cdots & & & 0 & 0 \\
\hline
-\mX_k & \lambda \mI_{\rho_k} -\mY_k & -\mZ_k & 0 & \cdots & & 0 & 0 \\
* & * & \lambda \mI_{\rho_{k-1}} -\mY_{k-1} & -\mZ_{k-1} & & & \vdots & 0 \\[8pt]
\vdots & \vdots & \ddots\hfill {} & \hfill\ddots & & \hfill \ddots & 0 & \vdots \\[8pt]
 & & & & * & \lambda \mI_{\rho_2} -\mY_2 & -\mZ_2 & 0 \\
* & * & & \cdots & & * & \lambda \mI_{\rho_1} -\mY_1 & \mZ_1
\end{array}\right] , 
\end{split}
\end{equation}
\new{where $\mI_\rho$ denotes the identity matrix in $\CNN{\rho}$, $\mA_k\in\CNN{\tau_k}$, $\mX_k\in\C^{\rho_k\times \tau_k}$, $\mY_j\in\CNN{\rho_j}$, $j=1,\ldots,k$; $\mZ_j\in\C^{\rho_j \times\rho_{j-1}}$, $j=2,\ldots,k$; $\mZ_1\in\C^{\rho_1 \times n}$ has full row rank $\rho_1$, and $*$ are matrices which are not computed in Algorithm~\ref{algorithm:staircase} below.
The constants $\tau_k$ and $\rho_j$, $j=1,\ldots,k$ indicate dimensions of submatrices; they are determined in Algorithm~\ref{algorithm:staircase} below and satisfy $\tau_k +\sum_{j=1}^k \rho_j =n$.} 

\begin{breakablealgorithm}
\caption{Staircase algorithm}
\label{algorithm:staircase}
\begin{algorithmic}[1]
\REQUIRE $c:=0$, $\mT:=\mI$, $\mA_0:=\mA$, $\mB_0:=\mB$, $j:=1$
\STATE Construct a unitary matrix~$\mU_j\in\Cnn$ such that $\displaystyle
\mU_j^H \mB_{j-1} =:
\left[\begin{array}{c} 0 \\ \hline \mZ_j \end{array}\right]
\begin{array}{c} \tau_j \\ \rho_j \end{array}$.
\label{algorithm:staircase:compression}
\STATE Extract $\rho_j$ and $\tau_j$ from Line~\ref{algorithm:staircase:compression}.
\WHILE{$\rho_j >0$ \algAnd $\tau_j >0$}
\STATE Transform and partition $\mA_{j-1}$ analogously such that
$\displaystyle
\mU_j^H \mA_{j-1} \mU_j
=:
\begin{array}{l}
\left[\begin{array}{c|c}
 \mA_j & \mB_j \\
 \hline
 \mX_j & \mY_j \\
\end{array}\right]
\begin{array}{c} \tau_j \\ \rho_j \end{array}
\\
\hspace{11pt} \tau_j \hspace{14pt} \rho_j
\end{array}
$.
\STATE Update
$\displaystyle
\mT
 := \mT
 \left[\begin{array}{c|c}
 \mU_j & 0 \\
 \hline
 0 & \mI_c
 \end{array}\right],
 \quad c:= c+\rho_j$.
\STATE $j:=j+1$
\STATE Construct a unitary matrix~$\mU_j\in\C^{\tau_{j-1}\times\tau_{j-1}}$ such that $\displaystyle
\mU_j^H \mB_{j-1} =:
\left[\begin{array}{c} 0 \\ \hline \mZ_j \end{array}\right]
\begin{array}{c} \tau_j \\ \rho_j \end{array}$.
\label{algorithm:staircase:compression:2}
\STATE Extract $\rho_j$ and $\tau_j$ from Line~\ref{algorithm:staircase:compression:2}.
\ENDWHILE
\IF{$\rho_j =0$}
\STATE $k:=j-1$ 
\hfill \COMMENT{Note that $k$ is the number of full rank ''stairs`` on termination of Algorithm~\ref{algorithm:staircase}.}
\ELSIF{$\tau_j =0$}
\STATE $k:=j$,\ $\mZ_k:=\mB_k$,\ $\mY_k:=\mA_k$ 
\ENDIF
\ENSURE Unitary matrix~$\mT$ such that~\eqref{staircase:VD} holds.
\end{algorithmic}
\end{breakablealgorithm}

\subsubsection*{Comparison of Algorithm~\ref{algorithm:1Pi} and  Algorithm~\ref{algorithm:staircase}}
\label{subsec:quadraticLyapunovfunctional}

To compare the Algorithms~\ref{algorithm:1Pi} and~\ref{algorithm:staircase}, we consider $\mA=\mJ$, $\mB=\mR$ such that $\mJ-\mR$ has finite HC-index~$\mHC$.
Algorithm~\ref{algorithm:staircase} constructs a unitary transformation to reduce a matrix pencil to the staircase form~\eqref{staircase:VD}, whereas Algorithm~\ref{algorithm:1Pi} uses these unitary transformations without computing the staircase form explicitly.
Both algorithms partition matrices.
Algorithm~\ref{algorithm:staircase} extracts submatrices, whereas Algorithm~\ref{algorithm:1Pi} uses projections, i.e., it is a coordinate free approach.
For a discussion of a coordinate free version of the SVD, see e.g.~\cite{OtPa15}.

To continue our comparison, we relate some of the matrices in  Algorithm~\ref{algorithm:staircase} and Algorithm~\ref{algorithm:1Pi}:
In Algorithm~\ref{algorithm:staircase}, the unitary transformation~$\mU_j$ of $\mB_{j-1}\in\C^{\tau_{j-1}\times\rho_{j-1}}$ is a complex matrix~$\mU_j\in\C^{\tau_{j-1}\times \tau_{j-1}}$ (with $\tau_0=n$, $\rho_0 =n$), where (a basis of) left-singular vectors corresponding to zero singular values of~$\mB_{j-1}$ are taken to be the first columns of~$\mU_j$.
In fact, we write
\begin{align} \label{Uj:partition}
\mU_j
&=
\left[\begin{array}{c|c} \mU_{j0} & \mU_{j1} \end{array}\right] \ \tau_{j-1} \\[-4pt]
&\hspace{24pt} \tau_j \hspace{18pt} \rho_j \nonumber
\end{align}
where the columns of $\mU_{j0}\in\C^{\tau_{j-1} \times \tau_j}$, resp. $\mU_{j1}\in\C^{\tau_{j-1} \times \rho_j}$, are orthogonal left-singular vectors of~$\mB_{j-1}$ corresponding to zero, resp. non-zero singular values of~$\mB_{j-1}$.
The left-singular vectors corresponding to a singular value~$0$ of~$\mB_{j-1}$ are elements of $\kernel\mB_{j-1}^H$.

Starting with $\mA_0 =\tA_0 =\mJ$ and $\mB_0 =\tB_0 =\mR$, for~$j=1,\ldots,\mHC$, we claim that the matrices such as~$\tA_j, \tB_j\in\Cnn$ in Algorithm~\ref{algorithm:1Pi} and the corresponding (sub)matrices~$\mA_j\in\C^{\tau_j\times \tau_j}, \mB_j\in\C^{\tau_j\times \rho_j}$ in Algorithm~\ref{algorithm:staircase} are related as
\begin{align}
\tA_j
&=(\mU_{10} \cdots\mU_{j0}) \mA_j (\mU_{10} \cdots\mU_{j0})^H \,,  \label{relation:A} \\
\tB_j
&=\begin{cases}
\mU_{10} \mB_1 \mU_{11}^H & \text{for } j=1 \,, \\
(\mU_{10} \cdots\mU_{j0}) \mB_j (\mU_{10} \cdots\mU_{j-1,0} \mU_{j1})^H & \text{for } j\geq 2 \,,
\end{cases} \label{relation:B} \\
\widetilde{\Pi}_j
&=\begin{cases}
\mU_{10} \mU_{10}^H & \text{if } j=1\,, \\
\mU_{11}\mU_{11}^H +(\mU_{10} \mU_{20})(\mU_{10} \mU_{20})^H & \text{if } j=2\,, \\
\mU_{11}\mU_{11}^H +\cdots +(\mU_{10}\cdots\mU_{j-2,0}\mU_{j-1,1})(\mU_{10}\cdots\mU_{j-2,0}\mU_{j-1,1})^H & \\
+(\mU_{10}\cdots\mU_{j,0})(\mU_{10}\cdots\mU_{j,0})^H & \text{if } j\geq 3\,,
\end{cases} \label{relation:tPi-U}
\end{align}
as well as
\begin{equation} \label{relation:Pi-U}
 \Pi_j
 =\widetilde{\Pi}_j \Pi_{j-1}
 =(\mU_{10}\cdots\mU_{j,0})(\mU_{10}\cdots\mU_{j,0})^H \,.
\end{equation}
The identities~\eqref{relation:A}--\eqref{relation:Pi-U} are proven iteratively:
Consider the Algorithms~\ref{algorithm:1Pi} and~\ref{algorithm:staircase} with initial data $\mA_0 =\tA_0 =\mJ$ and $\mB_0 =\tB_0 =\mR$.
Starting with~$j=1$, the projection matrices~$\widetilde\Pi_j,\Pi_j\in\Cnn$ and the unitary matrix~$U_j$ are computed from~$\tB_{j-1}$ and~$\mB_{j-1}$, respectively.
Using the partitioning of~$U_j$ as in~\eqref{Uj:partition} and the identities~\eqref{relation:A}--\eqref{relation:Pi-U} up to~$j-1$, the projection matrices~$\widetilde\Pi_j,\Pi_j$ and the unitary matrix~$\mU_j$ can be related as in~\eqref{relation:tPi-U}--\eqref{relation:Pi-U}.
Then the matrices~$\tA_j, \tB_j\in\Cnn$ in Algorithm~\ref{algorithm:1Pi} and the matrices~$\mA_j\in\C^{\tau_j\times \tau_j}, \mB_j\in\C^{\tau_j\times \rho_j}$ in Algorithm~\ref{algorithm:staircase} are computed and the identities~\eqref{relation:A}--\eqref{relation:B} are established.
This procedure is iterated until $j=\mHC$ is reached.

%

\section{Linear time-invariant DAE systems}
\label{sec:DAE}

\subsection{Canonical and condensed forms for semi-dissipative Hamiltonian DAEs}
\label{ssec:DAE:CF}

Two pairs of matrices~$(\mE_1,\mA_1), (\mE_2,\mA_2)\in\Cnn\times\Cnn$ are called \emph{equivalent} if there exist nonsingular matrices~$\mL,\mZ\in\Cnn$ such that $\mE_2 =\mL\mE_1\mZ$, $\mA_2=\mL\mA_1\mZ$.
If this is the case, we write $(\mE_1,\mA_1)\sim (\mE_2,\mA_2)$; see Definition~2.1 in~\cite{KunM06}.
This relation is an \emph{equivalence relation}, see Lemma~2.2 in~\cite{KunM06}, and it leads to the Kronecker canonical form~\cite{Gan59a}.
Let us denote by~$\mathcal J_k(\lambda_0)$ the standard upper triangular Jordan block of size $k\times k$ associated with the eigenvalue $\lambda_0$
and let $\mathcal L_k$ denote the standard right Kronecker block of size $k\times(k+1)$, i.e.,
\[
\mathcal J_k(\lambda_0)
=\begin{bmatrix} \lambda_0&1\\&\ddots&\ddots\\&&\ddots&1\\&&&\lambda_0 \end{bmatrix}
\quad\mbox{and}\quad
\mathcal L_k
=\lambda\begin{bmatrix} 1&0\\&\ddots&\ddots\\&&1&0 \end{bmatrix}
-\begin{bmatrix} 0&1\\&\ddots&\ddots\\&&0&1 \end{bmatrix} .
\]
\begin{theorem}[Kronecker canonical form, {Theorem 2.1 in \cite{MehMW18}}]\label{th:kcf}
Let $\mE,\mA\in {\mathbb C}^{n,m}$.
Then there exist nonsingular matrices $\mS \in {\mathbb C}^{n,n}$ and $\mT\in {\mathbb C}^{m,m}$ such that
\begin{equation}\label{kcf}
\mS (\lambda \mE-\mA)\mT
=\diag({\cal L}_{\eps_1},\ldots,{\cal L}_{\eps_p},
{\cal L}^\top_{\eta_1},\ldots,{\cal L}^\top_{\eta_q},
{\cal J}_{\rho_1}^{\lambda_1},\ldots,{\cal J}_{\rho_r}^{\lambda_r},
{\cal N}_{\sigma_1},\ldots,{\cal N}_{\sigma_s}),
\end{equation}
where $p,q,r,s\in\N_0$, $\eps_1,\dots,\eps_p,\eta_1,\dots,\eta_q,\rho_1,\dots,\rho_r,\sigma_1,\dots,\sigma_s\in\mathbb N$ and
$\lambda_1,\dots,\lambda_r\in\mathbb C$, as well as ${\cal J}_{\rho_i}^{\lambda_i}=\lambda \mI_{\rho_i}-\mathcal J_{\rho_i}(\lambda_i)$
for $i=1,\dots,r$ and $\mathcal N_{\sigma_j}=\lambda \mathcal J_{\sigma_j}(0)-\mI_{\sigma_j}$ for $j=1,\dots,s$.
This form is unique up to permutation of the blocks.
\end{theorem}
An eigenvalue~$\lambda$ is called semi-simple if the largest associated Jordan block~${\cal J}_\rho^\lambda$ has size one.
The sizes $\eta_j$ and $\eps_i$ of the rectangular blocks
are called the \emph{left and right minimal indices} of $\lambda \mE-\mA$, respectively.
The matrix pencil $\lambda \mE-\mA$ for $\mE,\mA \in \mathbb C^{n,m}$ is called \emph{regular} if $n=m$ and $\operatorname{det}(\lambda_0 \mE-\mA)\neq 0$ for some $\lambda_0 \in \mathbb C$, otherwise it is called \emph{singular}.
The values $\lambda_1,\dots,\lambda_r\in\mathbb C$ are called the finite eigenvalues of $\lambda \mE-\mA$.
If $s>0$, then $\lambda_0=\infty$ is said to be an eigenvalue of $\lambda \mE-\mA$.
The size of the largest block ${\cal N}_{\sigma_j}$ is
called the \emph{DAE-index (or Kronecker index)}~$\nu$ of the pencil~$\lambda \mE-\mA$, where, by convention,  $\nu=0$ if $\mE$ is invertible.

The spectral properties of a pencil~$\lambda \mE-\mA$ with semi-dissipative matrix~$\mA$ are characterized in the following theorem, which is a special case of results in~\cite{MehMW18,MehMW20}.
%
\begin{theorem}\label{thm:singind}
Let $\mE,\mJ,\mR\in\Cnn$ satisfying $\mE=\mE^H\geq 0$, $\mR=\mR^H\geq 0$ and $\mJ=-\mJ^H$.
Consider the matrix pencil~$P(\lambda)=\lambda \mE-(\mJ-\mR)$.
\begin{enumerate}[(E1)]
\item \label{lambda:non-pos}
If $\lambda_0\in\mathbb C$ is an eigenvalue of $P(\lambda)$ then $\Re(\lambda_0)\leq 0$.
\item \label{lambda:pure-imag}
If $\omega\in\mathbb R$ and $\lambda_0=i\omega$ is an eigenvalue of $P(\lambda)$, then $\lambda_0$ is semi-simple.
Moreover, if the columns of $\mV\in\mathbb C^{m,k}$ form a basis of a regular deflating subspace of $P(\lambda)$ associated with $\lambda_0$, then $\mR \mV=0$.
\item
The DAE-index of $P(\lambda)$ is at most two.
\item
All right minimal indices of $P(\lambda)$ are at most one (if there are any).
\item
All left minimal indices of $P(\lambda)$ are zero (if there are any).
\end{enumerate}
\end{theorem}
%
%
We will employ Theorem~\ref{thm:singind} for the case of regular pencils, i.e., when there are no left and right minimal indices.
In the special case~$\mE=\mI$, Theorem~\ref{thm:singind}~\ref{lambda:pure-imag} implies Proposition~\ref{prop:border}.

Employing  Lemma~\ref{lem:SF} we obtain the following staircase form.
\begin{lemma}[Staircase form for triple $(\mE,\mJ,\mR)$] \label{lem:tSF}
Let $\mE,\mJ,\mR\in\Cnn$ satisfy $\mE=\mE^H\geq 0$, $\mR=\mR^H\geq 0$ and $\mJ=-\mJ^H$.
Then there exists a unitary matrix $\mP\in\Cnn$, such that $\widecheck \mE :=\mP\ \mE\ \mP^H$, $\widecheck \mJ :=\mP\ \mJ\ \mP^H$ and $\widecheck \mR :=\mP\ \mR\ \mP^H$ satisfy
\begin{equation}\label{staircase:EJR}
\begin{split}
\widecheck \mE 
&=:\begin{bmatrix}
\mE_{1,1} & \mE_{2,1}^H & 0 & 0 & 0 \\
\mE_{2,1} & \mE_{2,2} & 0 & 0 & 0 \\
0  & 0 & 0 & 0 & 0 \\
0  & 0 & 0 & 0 & 0 \\
0  & 0 & 0 & 0 & 0
\end{bmatrix},
\\
\widecheck \mJ 
&=:\begin{bmatrix}
\mJ_{1,1} & -\mJ_{2,1}^H & -\mJ_{3,1}^H & -\mJ_{4,1}^H & 0\\
\mJ_{2,1} & \mJ_{2,2} & -\mJ_{3,2}^H& 0 & 0\\
\mJ_{3,1} & \mJ_{3,2} & \mJ_{3,3}&  0 & 0\\
\mJ_{4,1} & 0& 0 & 0 &0\\
 0 & 0 & 0 & 0 & 0
\end{bmatrix},
\qquad %
\widecheck \mR 
=:\begin{bmatrix}
\mR_{1,1} & \mR_{2,1}^H & \mR_{3,1}^H& 0 & 0 \\
\mR_{2,1} & \mR_{2,2} & \mR_{3,2}^H &0 & 0  \\
\mR_{3,1} & \mR_{3,2} & \mR_{3,3} &0 & 0 \\
0  & 0 & 0 & 0 & 0 \\
0  & 0 & 0 & 0 & 0
\end{bmatrix}.
\end{split}
\end{equation}
These three matrices are partitioned in the same way, with (square) diagonal block matrices of sizes $n_1,n_2,n_3,n_4=n_1,n_5\in\N_0$.
\new{If the block matrices $\mE_{1,1}$, $\mE_{2,2}$ (as well as $\mE_{2,1}$) are present, then the matrices $\mE_{1,1}$, $\mE_{2,2}$ (as well as $\begin{bmatrix} \mE_{1,1} & \mE_{2,1}^H \\ \mE_{2,1} & \mE_{2,2}  \end{bmatrix}$) are positive definite.
If the block matrices $\mJ_{4,1}$, $\mJ_{3,3} -\mR_{3,3}$ are present, then the matrices $\mJ_{4,1}$, $\mJ_{3,3}-\mR_{3,3}$ are invertible.} 
\end{lemma}
\begin{proof}
The proof is given as the following constructive Algorithm~\ref{algorithm:staircase:E_J-R}, which is similar to that of Lemma~\ref{lem:blocksvd}.
\begin{breakablealgorithm}
\caption{Staircase Algorithm for triple~$(\mE,\mJ,\mR)$}
\label{algorithm:staircase:E_J-R}
\begin{algorithmic}[1]
 ----------- \textit{Step 1} -----------
\STATE Perform a spectral decomposition of~$\mE$ such that
\[
\mE =\mU_E \begin{bmatrix} \tE_{1,1} & 0 \\ 0 & 0 \end{bmatrix} \mU_E^H,
\]
with $\mU_E\in\Cnn$ unitary, $\tE_{1,1}\in \mathbb C^{\tilde n_1,  \tilde n_1}$ positive definite or $\tilde n_1=0$.
\STATE Set $\mP := \mU^H_E$,
\begin{align*}
\tJ
&:= \mU^H_E\ \mJ\ \mU_E
=\begin{bmatrix}
 \tJ_{1,1} & -\tJ_{2,1}^H \\
 \tJ_{2,1} & \tJ_{2,2}
 \end{bmatrix},
& %
\tR
&:= \mU^H_E\ \mR\ \mU_E
=\begin{bmatrix}
 \tR_{1,1} & \tR_{2,1}^H \\
 \tR_{2,1} & \tR_{2,2} \end{bmatrix},
& %
\tE
&:=\mU^H_E\ \mE\ \mU_E.
\end{align*}
\newline ----------- \textit{Step 2} -----------
\IF{$\tilde n_1<n$}
\STATE Apply Lemma~\ref{lem:blocksvd} to $\tJ_{2,2}-\tR_{2,2}\in\C^{(n-\tilde n_1)\times (n-\tilde n_1)}$ such that
\[
\mP_{2,2}\ (\tJ_{2,2}-\tR_{2,2})\ \mP_{2,2}^H
=\begin{bmatrix} \tSigma_{2,2} & 0 \\ 0 & 0\end{bmatrix},
\]
with $\tSigma_{2,2}\in \mathbb C^{\tilde n_2,\tilde n_2}$ invertible or $\tilde n_2=0$.
\ENDIF
\STATE Set
\[
\mP_2:=\begin{bmatrix} \mI& 0 \\ 0 & \mP_{2,2} \end{bmatrix} \in\Cnn,\qquad \mP:=\mP_2 \mP.
\]
\STATE Set $\tE:=\mP_2\ \tE\ \mP_2^H$,
\[
\def\arraystretch{1.4}
\tJ := \mP_2\ \tJ\ \mP_2^H
=: \left[ \begin{array}{c|cc}
 \tJ_{1,1} & -\tJ_{2,1}^H & -\tJ_{3,1}^H  \\
 \hline 
 \tJ_{2,1} & \tJ_{2,2} & 0 \\
 \tJ_{3,1}& 0 & 0
\end{array}\right],\qquad
\tR := \mP_2\ \tR\ \mP_2^H
=: \left[ \begin{array}{c|cc}
 \tR_{1,1} & \tR_{2,1}^H & 0 \\
 \hline 
 \tR_{2,1} & \tR_{2,2} & 0 \\
 0 & 0 & 0
\end{array}\right],\quad
\]
with $\tJ_{2,2}-\tR_{2,2}=\tSigma_{2,2}$.
(The lines indicate the partitioning of the block matrices~$\tJ$ and~$\tR$ in the previous step.
Note that the positive semi-definiteness of the Hermitian matrix~$\mR$ implies the~$0$ structure in~$\tR$.)
\newline ----------- \textit{Step 3} -----------
\STATE Define $\tilde n_3 :=n -\tilde n_1 -\tilde n_2$.
\IF{$\tilde n_3>0$}
\STATE Perform an SVD of $\tJ_{3,1}$ such that
\[
\tJ_{3,1} = \mU_{3,1} \begin{bmatrix} \tSigma_{3,1} & 0\\ 0 & 0 \end{bmatrix} \mV^H_{3,1} \in \C^{\tilde n_3 \times \tilde n_1}\,,
\]
\hspace{12pt} with~$\tSigma_{3,1}\in\R^{n_1 \times n_1}$ nonsingular diagonal or $n_1 =0$.
\ENDIF
\STATE Set
\[
\mP_3 :=\begin{bmatrix}
 \mV_{3,1}^H & & \\
  & \mI & \\
  & & \mU_{3,1}^H
\end{bmatrix} \in\Cnn,\qquad
\mP:= \mP_3 \mP.
\]
\STATE Set $\widecheck \mE
:= \mP_3\ \tE\ \mP_3^H$, $\widecheck\mJ
:= \mP_3\ \tJ\ \mP_3^H$, $\widecheck\mR
:= \mP_3\ \tR\ \mP_3^H$ such that
\begin{align*}
\widecheck \mE
&=: \left[ \begin{array}{cc|c|cc}
 \mE_{1,1} &    \mE_{2,1}^H    &  0 &0 & 0 \\
 \mE_{2,1} & \mE_{2,2} &  0    &   0      & 0 \\
\hline
 0  & 0  & 0  & 0  & 0 \\
\hline
 0  & 0  & 0  & 0  & 0 \\
 0  & 0  & 0  & 0  & 0
\end{array}\right],
\\
\widecheck\mJ
&=:\left[ \begin{array}{cc|c|cc}
\mJ_{1,1}  &  -\mJ_{2,1}^H  &  -\mJ_{3,1}^H  &  -\mJ_{4,1}^H  &  0 \\
\mJ_{2,1}  &  \mJ_{2,2}  &  -\mJ_{3,2}^H  &  0  &  0 \\
\hline
\mJ_{3,1}  &  \mJ_{3,2}  &  \mJ_{3,3}  &  0  &  0 \\
\hline
\mJ_{4,1}  &  0  &  0  &  0  &  0 \\
0  &  0  &  0  &  0  &  0
\end{array}\right],
\qquad %
\widecheck\mR
=:\left[ \begin{array}{cc|c|cc}
\mR_{1,1}  &  \mR_{2,1}^H  &  \mR_{3,1}^H &  0  &  0 \\
\mR_{2,1}  &  \mR_{2,2}  &  \mR_{3,2}^H  &  0  &  0 \\
\hline
\mR_{3,1}  &  \mR_{3,2}  &  \mR_{3,3}    &  0  &  0 \\
\hline
0  &  0  &  0  &  0  &  0 \\
0  &  0  &  0  &  0  &  0
\end{array}\right],
\end{align*}
which are of the desired form with $n_2:=\tilde n_1 -n_1$, $n_3 :=\tilde n_2$, $n_4:=n_1$, $n_5:=\tilde n_3 -n_4$.
The matrices~$\mJ_{4,1}:=\tSigma_{3,1}$ and $\mJ_{3,3}-\mR_{3,3} =\tJ_{2,2}-\tR_{2,2} =\tSigma_{2,2}$ are invertible.
\end{algorithmic}
\end{breakablealgorithm}
\end{proof}

The pencil~$\lambda \mE-(\mJ-\mR)$ is associated to the DAE~\eqref{DAE:EA} with~$\mA=\mJ-\mR$, which can be transformed into a DAE in staircase form,
\begin{equation}\label{DAE:EA:hat:0}
\widecheck\mE\dot y =(\widecheck\mJ -\widecheck\mR)y\ ,
\qquad\text{with } y:=\mP x\ .
\end{equation}
Denoting $\widecheck\mA :=\widecheck\mJ -\widecheck\mR$, we shall prove next that $(\widecheck\mE,\widecheck\mA)$ is equivalent to some~$(\widehat\mE,\widehat\mA)$ whose pencil~$\lambda\widehat\mE -\widehat\mA$ is \emph{almost in Kronecker canonical form}, a terminology defined by~\eqref{kcf:almost} below.
Since $(\mE,\mA)\sim (\widecheck\mE,\widecheck\mA)$, also $(\mE,\mA)\sim (\widehat\mE,\widehat\mA)$.

\begin{lemma}\label{lem:reduced}
Consider a semi-dissipative Hamiltonian DAE in staircase form~\eqref{DAE:EA:hat:0} with~$\widecheck\mA:=\widecheck\mJ-\widecheck\mR$ using~\eqref{staircase:EJR}.
Then there exist nonsingular matrices~$\mL,\mZ$ such that
\begin{align}\label{kcf:almost}
\widehat\mE
&:=\mL\ \widecheck\mE\ \mZ =:
\begin{bmatrix}
\widehat\mE_{1,1} & 0 & 0 & 0 & 0 \\
0 & \widehat\mE_{2,2} & 0 & 0 & 0 \\
0 & 0 & 0 & 0 & 0 \\
0 & 0 & 0 & 0 & 0 \\
0 & 0 & 0 & 0 & 0
\end{bmatrix} ,
& %
\widehat\mA
&:=\mL\ \widecheck\mA\ \mZ
=:\begin{bmatrix}
0 & 0 & 0& \mI & 0 \\
0  & \widehat\mA_{2,2} & 0 & 0 & 0 \\
0 & 0 & \mI& 0 & 0 \\
-\mI & 0 & 0 & 0 & 0 \\
0 & 0 & 0 & 0 & 0
\end{bmatrix} .
\end{align}
The two matrices are partitioned in the same way, with (square) diagonal block matrices of sizes $n_1,n_2,n_3$, $n_4=n_1,n_5\in\N_0$.
If the matrices~$\widehat\mE_{1,1}$ and~$\widehat\mE_{2,2}$ are present, then they are Hermitian positive definite.
For example, if $n_1>0$, $n_2>0$, and $n_3>0$, then $\widehat\mE_{1,1} =\mE_{1,1} -\mE_{2,1}^H\mE_{2,2}^{-1}\mE_{2,1}$, $\widehat\mE_{2,2} =\mE_{2,2}$, and $\widehat\mA_{2,2} =\widecheck\mA_{2,2} -\widecheck\mA_{2,3}\widecheck\mA_{3,3}^{-1}\widecheck\mA_{3,2}$.
\end{lemma}
We call $\lambda\widehat\mE-\widehat\mA$ with $\widehat\mE, \widehat\mA$ given in~\eqref{kcf:almost} the \emph{almost Kronecker canonical form} of $\lambda\widecheck\mE-\widecheck\mA$, since the $\mathcal{L}$ and $\mathcal{N}$ blocks are already in the form given in~\eqref{kcf}.
\begin{proof}
\new{We discuss only the case that the block matrices~$\mE_{1,1}$, $\mE_{2,1}$ and $\mE_{2,2}$ are present in~\eqref{staircase:EJR}.
For all other cases the proof is similar, but the following block Gauss elimination is not needed.
If the block matrix~$\mE_{2,2}$ is present, then (due to Lemma~\ref{lem:tSF}) $\mE_{2,2}$ is positive definite.} 
To achieve the block diagonal structure of~$\widehat\mE$, perform block Gauss eliminations on the upper right diagonal block of matrix~$\widecheck\mE$ such that
\[
\begin{bmatrix}
\widehat\mE_{1,1} & 0 \\
0 &\widehat\mE_{2,2}
\end{bmatrix}
:=
\begin{bmatrix}
\mE_{1,1} -\mE_{2,1}^H\mE_{2,2}^{-1}\mE_{2,1} & 0 \\
0 & \mE_{2,2}
\end{bmatrix}
=
\underbrace{\begin{bmatrix} \mI & -\mE_{2,1}^H \mE_{2,2}^{-1} \\ 0 & \mI \end{bmatrix}}_{=:\mL_1}
\begin{bmatrix} \mE_{1,1} & \mE_{2,1}^H \\ \mE_{2,1} & \mE_{2,2} \end{bmatrix}
\underbrace{\begin{bmatrix} \mI & 0 \\ -\mE_{2,2}^{-1} \mE_{2,1} & \mI \end{bmatrix}}_{=:\mZ_1}
\ .
\]
The matrix $\widehat\mE_{1,1} =\mE_{1,1} -\mE_{2,1}^H\mE_{2,2}^{-1}\mE_{2,1}$ is a Schur complement of a positive definite Hermitian matrix, hence, it is again positive definite, see~\cite{JoSm05,JoSm06}.
Noting that $\mL_1 =\mZ_1^H$, the two matrices are congruent, hence, the Hermitian positive definiteness of~$\diag(\widehat\mE_{1,1},\widehat\mE_{2,2})$ follows from the Hermitian positive definiteness of the original matrix.

Defining the block diagonal matrices $\mL_\mE :=\diag(\mL_1,\mI)$ and $\mZ_\mE :=\diag(\mZ_1,\mI)$ yields $\breve\mE :=\mL_\mE\ \widecheck\mE\ \mZ_\mE =\diag(\widehat\mE_{1,1},\widehat\mE_{2,2},0,0,0)$ and
\begin{align*}
%
\breve\mA
:=\mL_\mE\ \widecheck\mA\ \mZ_\mE
=
&\begin{bmatrix}
\breve\mA_{1,1} & \widecheck\mA_{1,2} -\mE_{2,1}^H \mE_{2,2}^{-1}\widecheck\mA_{2,2} & \widecheck\mA_{1,3} -\mE_{2,1}^H \mE_{2,2}^{-1}\widecheck\mA_{2,3} & \widecheck\mA_{1,4} & 0 \\
\widecheck\mA_{2,1} -\widecheck\mA_{2,2} \mE_{2,2}^{-1} \mE_{2,1} & \widecheck\mA_{2,2} & \widecheck\mA_{2,3} & 0 & 0 \\
\widecheck\mA_{3,1} -\widecheck\mA_{3,2}\mE_{2,2}^{-1} \mE_{2,1} & \widecheck\mA_{3,2} & \widecheck\mA_{3,3} & 0 & 0 \\
\widecheck\mA_{4,1} & 0& 0 & 0 & 0 \\
0 & 0 & 0 & 0 & 0
\end{bmatrix} %
\\
=:
&\begin{bmatrix}
\breve\mA_{1,1} & \breve\mA_{1,2} & \breve\mA_{1,3} & \breve\mA_{1,4} & 0 \\
\breve\mA_{2,1} & \breve\mA_{2,2} & \breve\mA_{2,3} & 0 & 0 \\
\breve\mA_{3,1} & \breve\mA_{3,2} & \breve\mA_{3,3} & 0 & 0 \\
\breve\mA_{4,1} & 0& 0 & 0 & 0  \\
0 & 0 & 0 & 0 & 0
\end{bmatrix}
\ ,
\end{align*}
where $\breve\mA_{1,1} :=\widecheck\mA_{1,1} -\widecheck\mA_{1,2} \mE_{2,2}^{-1} \mE_{2,1} -\mE_{2,1}^H \mE_{2,2}^{-1} (\widecheck\mA_{2,1} -\widecheck\mA_{2,2} \mE_{2,2}^{-1} \mE_{2,1})$.
Using the special structure of the matrices~$\breve\mE$ and~$\breve\mA$, the next transformation matrices are
\[
\arraycolsep=4pt
\def\arraystretch{1.2}
\mL_\mA
:=\left[\begin{array}{ccccc}
\mI & 0 & -\breve\mA_{1,3} \breve\mA_{3,3}^{-1} & (-\breve\mA_{1,1}+ \breve\mA_{1,3} \breve\mA_{3,3}^{-1} \breve\mA_{3,1}) \breve\mA_{4,1}^{-1} & 0 \\
0 & \mI & -\breve\mA_{2,3} \breve\mA_{3,3}^{-1} & (-\breve\mA_{2,1}+ \breve\mA_{2,3} \breve\mA_{3,3}^{-1} \breve\mA_{3,1}) \breve\mA_{4,1}^{-1} & 0 \\
 0 & 0 & \breve\mA_{3,3}^{-1} & -\breve\mA_{3,3}^{-1} \breve\mA_{3,1} \breve\mA_{4,1}^{-1} & 0 \\
0 & 0 & 0 & -\breve\mA_{4,1}^{-1} & 0 \\
0 & 0 & 0 & 0 & \mI
\end{array}\right]
\]
and
\[
\mZ_\mA
:=\begin{bmatrix}
\mI & 0 & 0 & 0 & 0 \\
0 & \mI & 0 & 0 & 0 \\
0 & -\breve\mA_{3,3}^{-1} \breve\mA_{3,2}& \mI & 0 & 0 \\
0 & \breve\mA_{1,4}^{-1} (-\breve\mA_{1,2}+\breve\mA_{1,3}\breve\mA_{3,3}^{-1} \breve\mA_{3,2})  & 0 & \breve\mA_{1,4}^{-1} & 0 \\
0 & 0 & 0 & 0 & \mI
\end{bmatrix} \ .
\]
Altogether, the transformation matrices $\mL:=\mL_\mA \mL_\mE$ and $\mZ:=\mZ_\mE \mZ_\mA$ yield the equivalence of $(\widehat\mE,\widehat\mA)\sim (\widehat\mE,\widehat\mA)$ as stated in~\eqref{kcf:almost}.
\end{proof}
Note \new{that the block-sizes in both Lemma~\ref{lem:tSF} and~\ref{lem:reduced} are the same, and}
that we retain the dissipative structure of the dynamic equations associated with~$\mE_{2,2}$, $\mA_{2,2}$ despite the fact that we are performing nonunitary equivalence transformations in Lemma~\ref{lem:reduced}.
But these transformations only affect the algebraic equations:
%
%
\begin{corollary}\label{cor:index}
Let $\mE,\mJ,\mR\in\Cnn$ satisfy $\mE=\mE^H\geq 0$ with~$\mE\ne 0$, $\mR=\mR^H\geq 0$ and $\mJ=-\mJ^H$.
Consider the pencil~$\lambda \mE-(\mJ-\mR)$, its unitarily congruent pencil~$\lambda \widecheck \mE -(\widecheck \mJ-\widecheck \mR)$ in staircase form~\eqref{staircase:EJR} for some unitary matrix~$\mP\in\Cnn$, and its equivalent pencil~$\lambda\widehat\mE -\widehat\mA$ in almost Kronecker canoncial form~\eqref{kcf:almost} for some invertible matrices $\mL,\mZ\in\Cnn$ with parameters~$n_1, n_2, n_3, n_4, n_5\in\N_0$.
\begin{enumerate}[(P1)]
\item \label{P:n5}
The pencil~$\lambda \mE-(\mJ-\mR)$ is regular if and only if~$n_5=0$.
\item \label{P:DAE-index}
The DAE-index~$\nu$ of a regular pencil~$\lambda \mE-(\mJ-\mR)$ satisfies
\begin{equation}
\nu
=\begin{cases}
2 &\text{if and only if~$n_1=n_4>0$,} \\
1 &\text{if and only if~$n_1=n_4=0$ and $n_3>0$,} \\
0 &\text{if and only if~$n_1=n_4=0$ and $n_3=0$.}
\end{cases}
\end{equation}
\end{enumerate}

The pencil~$\lambda \mE-(\mJ-\mR)$ is associated to the DAE~\eqref{DAE:EA} with~$\mA=\mJ-\mR$, which can be transformed into a DAE in staircase form,
\begin{equation}\label{DAE:EA:UnitarilyCongruent}
\widecheck\mE\dot y =(\widecheck\mJ -\widecheck\mR)y\ ,
\qquad\text{with } y:=\mP x\ ,
\end{equation}
and furthermore into a DAE in almost Kronecker canonical form
\begin{equation}\label{DAE:EA:Equivalent}
\widehat\mE\dot z =\widehat\mA z\ ,
\qquad\text{with } z:=\mZ^{-1} y\ .
\end{equation}
\begin{enumerate}[resume*]
\item \label{DAE:EA:condensedForms}
If a regular pencil~$\lambda \mE-(\mJ-\mR)$ has DAE-index two, then we can eliminate the first and fourth equation in~\eqref{DAE:EA:UnitarilyCongruent} and~\eqref{DAE:EA:Equivalent}.
For example, if $n_2>0$ and $n_3>0$ then we obtain
\begin{align}\label{DAE:EA:kcf:almost}
z_1 =0\ , \qquad
\widehat\mE_{2,2} \dot z_2 =\widehat\mA_{2,2} z_2\ , \qquad
z_3 =0\ , \qquad
z_4 =0\ ,
\end{align}
which implies via $y=\mZ z$ that
\begin{equation}\label{DAE:EA:UC:coordinates}
\begin{split}
 y_1 &=0\ , \qquad
 y_2  =z_2\ , \qquad
 y_3  =-\widecheck\mA_{3,3}^{-1} \widecheck\mA_{3,2} y_2\ ,
\\
 y_4
&=\mJ_{4,1}^{-H}\left ( (-\mJ_{2,1}^H-\mR_{2,1}^H) y_2 +(-\mJ_{3,1}^H -\mR_{3,1}^H) y_3 -\mE_{2,1}^H\dot y_2\right ) \ ,
\end{split}
\end{equation}
leading to restrictions in the initial values.
\item \label{EJR:ODE}
For systems~\eqref{DAE:EA:UnitarilyCongruent} and~\eqref{DAE:EA:Equivalent} with~$n_5=0$, the underlying implicit ODE systems are given by the system in~$y_2$ and~$z_2$, respectively, that are obtained by eliminating all other variables.
For example, if $n_2>0$ and $n_3>0$, then this yields systems of the form
\begin{equation}\label{underlyingODE:z}
\widehat\mE_{2,2} \dot z_2
=\widehat\mA_{2,2} z_2,
\end{equation}
or equivalently
\begin{equation}\label{underlyingODE:y}
\mE_{2,2} \dot y_2
=\widehat \mA_{2,2} y_2
=(\widehat \mJ_{2,2}-\widehat \mR_{2,2}) y_2,
\end{equation}
with $\mE_{2,2}=\widehat\mE_{2,2}$ Hermitian positive definite and $\widehat \mA_{2,2}$ semi-dissipative.
Here
\[
 \widehat\mJ_{2,2} :=(\widehat\mA_{2,2})_S \ , \qquad
 \widehat\mR_{2,2} :=-(\widehat\mA_{2,2})_H \ , \quad
 \text{where } \widehat\mA_{2,2} =\widecheck\mA_{2,2} -\widecheck\mA_{2,3}\widecheck\mA_{3,3}^{-1}\widecheck\mA_{3,2} \ .
\]
%
\item
The finite eigenvalues of the system are the eigenvalues of the matrix pencil $\lambda \mE_{2,2} -\widehat \mA_{2,2}$.
They are in the closed left half plane and the eigenvalues on the imaginary axis are semi-simple.
\end{enumerate}
\end{corollary}
%
%
\begin{proof}
The proof follows directly from the analysis of the DAE, see~\cite{KunM06}.
In part~\ref{EJR:ODE}, $\mE_{2,2}$ is positive definite, since it is a principal submatrix of a positive definite matrix, whereas~$\widehat \mA_{2,2}$ is semi-dissipative, since Schur complements of semi-dissipative matrices again have this property, see~\cite{JoSm05,JoSm06}. 
\end{proof}

\begin{remark}\label{rem:rank}
It should be noted that the construction of the staircase form~\eqref{staircase:EJR} requires three consecutive rank decisions, which even though they are performed via unitary transformations, are critical from the point of view of perturbation theory.
It is an open problem to perform a detailed perturbation analysis how the different rank decisions affect each other.
In the context of the perturbation theory for general DAEs and staircase forms, it has been observed that it is better to consider a smaller rank if the rank decision is doubtful,~\cite{ByeGM97,KunM06}.
In the recent paper~\cite{MehMW20} a characterization of the smallest perturbation that makes the system singular or of higher index or unstable (which is the same) has been characterized via a simple optimization problem, so that these distances are accessible.
\end{remark}

\subsection{Negative hypocoercive pencils for semi-dissipative Hamiltonian DAEs}
\label{sec:DHDAE}

Many of the results on (negative) hypocoercive matrices from the previous section can be extended to linear semi-dissipative Hamiltonian DAEs of the form~\eqref{DAE:EA}.
These results will be based on the staircase form~\eqref{staircase:EJR} for the triple~$(\mE,\mJ,\mR)$ and the underlying implicitly defined semi-dissipative Hamiltonian ODE~\eqref{underlyingODE:y}.
It follows directly from the staircase form that the initial values have to be chosen in a consistent way.
If the system~\eqref{DAE:EA} is regular and is transformed to a DAE in staircase form~\eqref{DAE:EA:UnitarilyCongruent} using~\eqref{staircase:EJR}, then the transformed initial conditions for~$y_2$ can be chosen arbitrarily while the initial conditions for the other components $y_1, y_3, y_4$ are then already fixed, see e.g.~\ref{DAE:EA:condensedForms}.
To extend the definition of (negative) hypocoercivity to the DAE case, we therefore have to assume that the system is regular and that the initial conditions are consistent.
\begin{remark}\label{rem:IC:consistency}
Note that, for homogeneous linear time-invariant DAE systems, the set of consistent initial conditions is a linear subspace, see Corollary 2.30 in~\cite{KunM06}.
\end{remark}
One has to be careful if the system is of index two because then arbitrary small general perturbations can make the system unstable, however, if the perturbation stays in the class of semi-dissipative Hamiltonian DAEs then the perturbed system remains stable.
Consider the following example modified from~\cite{DuLM13}.

\begin{eexample} \label{ex:DAE:index2}
The homogeneous linear time-invariant DAE system $\dot x= y$, $0=-x-\varepsilon y$, with $\varepsilon \in\R$, is semi-dissipative if~$\varepsilon \geq 0$.
If~$\varepsilon >0$ then the DAE has index one and has the solution $x(t)=e^{-\varepsilon^{-1} t}x(0)$, $y(t)= -\varepsilon^{-1}e^{-\varepsilon^{-1} t}x(0)$.
With a consistent initial value $y(0)= -\varepsilon^{-1} x(0)$, the solution is asymptotically going to~$0$ for every initial condition~$x(0)$.
If $\varepsilon=0$ then the DAE has index two and the solution is $x=y=0$, now~$x(0)$ is restricted as well.
However, if $\varepsilon <0$ then the system is unstable.
\end{eexample}
\begin{definition}\label{def:hypoDAE}
A matrix pencil $\lambda \mE-\mA$ is called~\emph{negative hypocoercive} if the pencil is regular, of DAE-index at most two and the finite eigenvalues of the pencil $\lambda \mE-\mA$ have negative real part.
\end{definition}
We note that a regular pencil might not have \emph{any} finite eigenvalues, in which case the last condition would be void.
Due to Theorem~\ref{thm:singind}, a linear semi-dissipative Hamiltonian DAE system~\eqref{DAE:EA} with a regular pencil~$\lambda \mE-\mA$ only has finite eigenvalues with non-positive real part.
\begin{definition} \label{def:DAE:mHC}
Consider a linear semi-dissipative Hamiltonian DAE system~\eqref{DAE:EA} with a regular pencil~$\lambda \mE-\mA$ and the unitarily congruent DAE~\eqref{DAE:EA:UnitarilyCongruent} in staircase form~\eqref{staircase:EJR}.
If the underlying implicit ODE~\eqref{underlyingODE:y} is missing (present) then system~\eqref{DAE:EA} is said to exhibit \emph{(non-)trivial dynamics}.
In case of non-trivial dynamics, the~\emph{HC-index~$m_{HC}$} of~$\lambda \mE-\mA$ is defined as the HC-index of the system matrix~$(\mE_{2,2}^{1/2})^{-1} \widehat{\mA}_{2,2} (\mE_{2,2}^{1/2})^{-1}$ of~\eqref{underlyingODE:y}, otherwise it is defined as~$0$.
\end{definition}
\begin{remark} \label{remark:DAE:HCI}
\begin{enumerate}[(i)]
\item
The definition of the HC-index for semi-dissipative Hamiltonian DAE systems~\eqref{DAE:EA} with non-trivial dynamics can also be based on the equivalent DAE~\eqref{DAE:EA:Equivalent} in almost Kronecker form~\eqref{kcf:almost}, since the underlying implicit ODE systems~\eqref{underlyingODE:y} and~\eqref{underlyingODE:z} are identical:
\[
 y_2 =z_2 \ , \qquad
 \mE_{2,2} =\widehat\mE_{2,2} \ , \qquad
 \widehat \mA_{2,2}  =\widehat{\mJ}_{2,2} -\widehat{\mR}_{2,2} \ ,
\]
due to Lemma~\ref{lem:reduced} and Corollary~\ref{cor:index}.
\item
Already in the case of semi-dissipative Hamiltonian DAEs~\eqref{DAE:EA} with positive definite Hermitian matrix~$\mE$, the definition of the HC-index has to take into account~$\mA$ \emph{and}~$\mE$, since it is natural to base its definition on the HC-index of an equivalent ODE (see e.g.~\eqref{ODE:A-tilde} or~\eqref{ODE:E-1A}). 
Note that, for $\mE\in\PDHn$ fixed and general semi-dissipative~$\mA$, considering all equivalence transformations of the matrix pencil such that $(\mE,\mA)\sim (\mI,\widetilde\mA)$, congruence transformations like $\tA =\msqrtEinv\mA\msqrtEinv$ preserve the semi-dissipativity of~$\mA$.
Hence, Definition~\ref{def:DAE:mHC} is based on~$(\mE_{2,2}^{1/2})^{-1} \widehat{\mA}_{2,2} (\mE_{2,2}^{1/2})^{-1}$.
\item \label{remark:DAE:HCI:semi-norm:tX}
Continuing with the case of a semi-dissipative Hamiltonian DAE~\eqref{DAE:EA} with $\mE\in\PDHn$, the simple energy estimate~\eqref{energy:estimate} may not allow to prove stability.
If we consider instead a squared weighted norm~$\|\cdot\|_{\tX}^2$ for some $\tX\in\PDHn$ then, for solutions~$x(t)$ of~\eqref{DAE:EA}, we deduce
\begin{equation} \label{energy:estimate:X}
\ddt \|x(t)\|_{\tX}^2
= \ip{x(t)}{(\mA^H \mE^{-1}\tX +\tX\mE^{-1}\mA)x(t)} \ .
\end{equation}
If a positive definite Hermitian matrix~$\tX$ satisfies the Lyapunov matrix inequality
\begin{equation} \label{LMI:tildeX}
 \mA^H \mE^{-1}\tX +\tX\mE^{-1}\mA
 \leq 0 \ ,
\end{equation}
then the weighted norm of solutions $\|x(t)\|_{\tX}$ would decay monotonically for $t\geq 0$.
For example, the choice $\tX=\mE$ allows to conclude
\begin{equation} \label{energy:estimate:E}
\ddt \|x(t)\|_{\mE}^2
= \ip{x(t)}{(\mA^H +\mA)x(t)}
\leq 0 \ ,
\end{equation}
due to the semi-dissipativity of~$\mA$.
Note that~\eqref{energy:estimate:E} also holds for $\mE=\mE^H$ only positive semi-definite; in this case $\|x\|_{\mE}$ is only a semi-norm.
\end{enumerate}
\end{remark}
The following proposition is the DAE-counterpart of Theorem~\ref{th:HC-decay}.
It says that the HC-index again characterizes the short time behavior of its solution propagator, but restricted to the dynamical subspace.
\begin{proposition}\label{prop:DAE+HC-decay}
Consider the semi-dissipative Hamiltonian DAE~\eqref{DAE:EA} with a regular, negative hypocoercive pencil~$\lambda\mE-\mA$, DAE-index at most two, non-trivial dynamics, and consistent initial condition~$x(0)$, see Remark~\ref{rem:IC:consistency}.
Then its (finite) HC-index is $\mHC\in\N_0$, if and only if
\[
\|S(t)\|_{\mE}
=
1 -c t^a +\bigO(t^{a+1})
\qquad
\text{for } t\to 0+ \ ,
\]
where $c>0$ and $a=2\mHC+1$, and the propagator (semi-)norm pertaining to the evolution of~\eqref{DAE:EA} reads
\[
\|S(t)\|_{\mE}
:=
\sup_{\stackrel{\|x(0)\|_{\mE}=1}{\text{for consistent } x(0)}}
\|x(t)\|_{\mE} \ ,
\qquad
t\geq 0 \ .
\]
\end{proposition}
\begin{proof}
In the following computation we use~\eqref{DAE:EA:UnitarilyCongruent}, \eqref{staircase:EJR}, \eqref{DAE:EA:UC:coordinates}:
\[
\|x\|_{\mE}^2
=x^H \mE x
=y^H \mP\mE\mP^H y
=y^H \widecheck{\mE} y
=\begin{bmatrix} y_1^H, y_2^H \end{bmatrix}
\begin{bmatrix} \mE_{1,1} & \mE_{2,1}^H \\ \mE_{2,1} & \mE_{2,2} \end{bmatrix}
\begin{bmatrix} y_1 \\ y_2 \end{bmatrix}
=y_2^H \mE_{2,2} y_2
=\|w_2\|^2 \ ,
\]
where $w_2:=\mE_{2,2}^{1/2} y_2 \in\C^{n_2}$ satisfies due to~\eqref{underlyingODE:y} the ODE
\begin{equation}\label{ODE:w2}
\dot{w}_2
=
(\mE_{2,2}^{1/2})^{-1} \widehat{\mA}_{2,2} (\mE_{2,2}^{1/2})^{-1} w_2 \ .
\end{equation}
With Theorem~\ref{th:HC-decay} we obtain
\[
\|S(t)\|_{\mE}
=
\sup_{\|w_2(0)\|=1} \|w_2(t)\|
=
1 -c t^a +\bigO(t^{a+1})
\qquad
\text{for } t\to 0+ \ .
\]
Here, $a=2\mHC+1$, and $\mHC$ is the HC-index of the system matrix in~\eqref{ODE:w2}, and by Definition~\ref{def:DAE:mHC} also of~$\lambda\mE-\mA$.
\end{proof}
The following example shows that the congruence transformation $\mA\mapsto (\mE^{1/2})^{-1} \mA (\mE^{1/2})^{-1}$ may change the HC-index.
Hence, the ODE~\eqref{ODE:A} and the DAE~\eqref{DAE:EA} with the same matrix~$\mA$ may have different HC-indices, even for $\mE=\mE^H>0$.
\begin{eexample}\label{ex:DAE+HCI}
Consider the DAE~\eqref{DAE:EA} with $\mA=\mJ-\mR$ given in Example~\ref{ex:HCmatrix+perturbation} and the following positive definite Hermitian matrix~$\mE$,
\begin{align}\label{matrix:E:4x4}
\mE
&:=\begin{bmatrix}
 1 & 0 & 0 & 0 \\
 0 & \sfrac{20}9 & 0 &\sfrac{-16}9 \\
 0 & 0 & 1 & 0 \\
 0 &\sfrac{-16}9 & 0 & \sfrac{20}9
\end{bmatrix}; \text{ hence }
 & %
\mE^{-1}
&=
\begin{bmatrix}
 1 & 0 & 0 & 0 \\
 0 & \sfrac54 & 0 & 1 \\
 0 & 0 & 1 & 0 \\
 0 & 1 & 0 & \sfrac54
\end{bmatrix},
 & %
\msqrtEinv
&=\begin{bmatrix}
 1 & 0 & 0 & 0 \\
 0 & 1 & 0 & \sfrac12 \\
 0 & 0 & 1 & 0 \\
 0 & \sfrac12 & 0 & 1
\end{bmatrix} .
\end{align}
The semi-dissipative matrix~$\mA$ has HC-index~2, see Example~\ref{ex:HCmatrix+perturbation}.
By contrast, $\tA =\msqrtEinv \mA \msqrtEinv =:\tJ-\tR$ has HC-index~1 since
\begin{align*}
\tJ
&=\begin{bmatrix}
 0 & 1 & 0 & \sfrac12 \\
-1 & 0 & 1 & 0 \\
 0 &-1 & 0 & \sfrac{-1}2 \\
\sfrac{-1}2 & 0 & \sfrac12 & 0
\end{bmatrix},
 & %
\tR
&=\begin{bmatrix}
 0 & 0 & 0 & 0 \\
 0 & \sfrac14 & 0 & \sfrac12 \\
 0 & 0 & 1 & 0 \\
 0 & \sfrac12 & 0 & 1
\end{bmatrix} \in\PSDHn ,
 & %
\tR +\tJ\tR\tJ^H
&=\begin{bmatrix}
 1 & 0 &-1 & 0 \\
 0 & \sfrac54 & 0 & 1 \\
-1 & 0 & 2 & 0 \\
 0 & 1 & 0 & \sfrac54
\end{bmatrix} \in\PDHn .
\end{align*}
\end{eexample}
The next example shows that for semi-dissipative Hamiltonian DAEs~\eqref{DAE:EA} with $\mE\in\PSDHn$, the squared weighted semi-norm~$\|\cdot\|_\mE^2$ only captures the behavior in the dynamic component~$y_2$ but not in components~$y_3, y_4$:
\begin{eexample} \label{ex:DAE:seminorm}
For $\varepsilon>0$, consider a linear semi-dissipative Hamiltonian DAE system~\eqref{DAE:EA} in staircase form~\eqref{DAE:EA:UnitarilyCongruent} with matrices
\begin{align*}
\widecheck\mE
&:=\begin{bmatrix}
  5 & 3 & 0 & 0 \\
  3 & 2 & 0 & 0 \\
  0 & 0 & 0 & 0 \\
  0 & 0 & 0 & 0
 \end{bmatrix} ,
& %
\widecheck\mJ
&:=\begin{bmatrix}
  0 & 1 & 1 & \varepsilon \\
 -1 & 0 & 1 & 0 \\
 -1 &-1 & 0 & 0 \\
 -\varepsilon & 0 & 0 & 0
 \end{bmatrix} ,
& %
\widecheck\mR
&:=\begin{bmatrix}
  0 & 0 & 0 & 0 \\
  0 &-2 & 0 & 0 \\
  0 & 0 &-1 & 0 \\
  0 & 0 & 0 & 0
 \end{bmatrix} ,
\end{align*}
such that $n_1=n_2=n_3=n_4=1$.
The associated pencil~$\lambda \mE-\mA$ is regular and exhibits nontrivial dynamics.
Due to Corollary~\ref{cor:index}~\ref{DAE:EA:condensedForms}--\ref{EJR:ODE}, we find that
\[
 y_1(t)=0 , \qquad
 \dot y_2 (t) =-y_2(t) , \qquad
 y_3(t) =-y_2(t) , \qquad
 y_4(t) =-\tfrac3{\varepsilon} y_2(t) \,.
\]
Thus, for given~$y_2(0)\in\R$, the solution of~\eqref{DAE:EA} is
\[
 y_1(t)=0 , \qquad
 y_2(t) = y_2(0)\ e^{-t} , \qquad
 y_3(t) =-y_2(t) =-y_2(0)\ e^{-t} , \qquad
 y_4(t) =-\tfrac3{\varepsilon}  y_2(0)\ e^{-t} ,
\]
such that $y_4(0) =-3 y_2(0)/{\varepsilon}$ can be arbitrarily large for sufficiently small $\varepsilon>0$.
In contrast, the squared weighted semi-norm of this solution satisfies $\|y(t)\|_\mE^2 =2 (y_2(0))^2 e^{-2t}$ for $t\geq 0$.
\end{eexample}
Note that for DAEs the concept of asymptotic stability is defined differently in the literature.
Often it is required that the DAE-index is at most~1, because otherwise there are hidden consistency conditions for the initial values and smoothness requirements for inhomogeneities.
Since we are only discussing linear homogeneous problems we allow DAE-index two.
\begin{corollary}\label{cor:hypodae}
If a semi-dissipative Hamiltonian DAE of the form~\eqref{DAE:EA} has a regular pencil $\lambda\mE-(\mJ-\mR)$ with DAE-index at most two, and non-trivial dynamics with a finite HC-index, then for every consistent initial condition the solution is asymptotically stable.
This property is retained for every sufficiently small perturbation of the system that stays within the class of semi-dissipative Hamiltonian DAEs, regardless whether the DAE-index changes, as long as the associated pencil stays regular.
\end{corollary}
\begin{proof}
By performing a congruence transformation, we may assume that the system is in almost Kronecker form~\eqref{DAE:EA:Equivalent} using~\eqref{kcf:almost}.
The regularity of the pencil ensures that~$n_5=0$ and that the DAE-index is at most two.
Moreover, the system~\eqref{DAE:EA:Equivalent} has the form~\eqref{DAE:EA:kcf:almost}, where some components may be missing.

Due to the assumption of non-trivial dynamics with finite HC-index, it follows that the solution~$z_2(t)$ of~\eqref{underlyingODE:z} is asymptotically going to zero for every initial value~$z_2(0)$.
Then the solution of the original system~$\mE\dot x=(\mJ-\mR)x$ is $x =\mP^H \mZ [0, z_2, 0, 0]^\top$, hence, it asymptotically goes to zero:
\begin{equation*} 
 \| x(t)\|^2
 =\| \mP^H \mZ z(t) \|^2
 =\| \mZ z(t) \|^2
 \leq \sigma_{\max}(\mZ) \|z_2(t)\|^2
 \leq \sigma_{\max}(\mZ) \|z_2(0)\|^2 e^{-2\mu t}
 \leq \kappa(\mZ) \|x(0)\|^2 e^{-2\mu t} ,
\end{equation*}
where $x(0)$ is a consistent initial value of the form $x(0) =\mP^H \mZ [0, z_2(0), 0, 0]^\top$, $\sigma_{\max}(\mZ)$ is the largest singular value of~$\mZ$, $\mu>0$ is some exponential decay rate capturing the asymptotic stability of~\eqref{underlyingODE:z}, and $\kappa(\mZ)=\|\mZ\| \|\mZ^{-1}\|$ is the condition number of~$\mZ$.


If the system is perturbed within the class of semi-dissipative Hamiltonian DAEs and the perturbed system is still regular, then the DAE-index may change between zero, one, or two and the set of consistent initial condition changes as well, but whenever they stay consistent, the asymptotic stability stays invariant.
\end{proof}

\begin{remark}\label{rem:singularity} Note that the solution~$x =\mP^H \mZ [0, z_2, 0, 0]^\top$ involves the inverses of the matrices $\mE_{2,2}$, $(\mJ_{3,3}-\mR_{3,3})$, and $\mJ_{4,1}$.
If these matrices are close to being singular, then this will lead to very large solution components in~$x$.
Thus, even though~$x$ asymptotically goes to~$0$, the decay may start from a very large level and so it may take a very long time until the solution is close to zero, even if the finite eigenvalues have large negative real parts, see Examples~\ref{ex:DAE:index2} and~\ref{ex:DAE:seminorm}.
\end{remark}

\begin{remark}\label{rem:daecontr} 
\new{The relationship between hypocoercivity and controllability for semi-dissipative \new{Hamiltonian} DAEs is analogous as in Lemma~\ref{lem:Equivalence}. Once the algebraic equations have been separated from the dynamic equations as in Lemma~\ref{lem:tSF}, one just applies Lemma~\ref{lem:Equivalence} to the dynamic part. 
Note, however, that Kalman matrix conditions like \ref{B:KRC} are not defined in the DAE case, one rather uses condition \ref{B:PBH:EVec} which then reads: No generalized eigenvector of~$\lambda \mE-\mJ$  associated with a finite eigenvalue lies in the kernel of~$\mR$, and 
condition \ref{B:PBH:EVal} which then reads:
$\rank [\lambda \mE-\mJ, \mR] =n$ for every $\lambda \in \C$, in particular for every finite eigenvalue~$\lambda$ of the matrix pencil~$\lambda \mE-\mJ$.
}
\end{remark}

\subsection{Negative hypocoercivity of DAEs and Lyapunov stability}\label{ssec:DAE:lyastab}

In the case of semi-dissipative Hamiltonian DAEs with singular~$\mE$, the situation is much more complex than for regular~$\mE$, since the classical relation between the existence of positive definite solutions to Lyapunov equations and stability of a system does not hold any longer, see~\cite{Sty02,Sty02a} for a detailed analysis.
The essential difference is that the solution of the Lyapunov equation need not be semi-definite, only the part associated with the dynamic part, and furthermore the right hand side has to be adapted.

Consider a linear DAE~$\mE \dot x=\mA x$ with square matrices~$\mE,\mA \in \Cnn$ and an associated generalized Lyapunov equation
\begin{equation}\label{eq:genlya}
\mE^H \mX \mA+ \mA^H \mX\mE=-\mE^H \mW \mE
\end{equation}
for some~$\mX,\mW\in\Cnn$.
This system has been studied e.g. in~\cite{Lew86,Mae98,ReReVo15,Sty02,Sty02a}.
The results in~\cite{Sty02,Sty02a} imply the following theorem.
\begin{theorem}\label{thm:genlya}
Consider a semi-dissipative Hamiltonian DAE~\eqref{DAE:EA} whose matrix pencil~$\lambda \mE-\mA$ is regular and has finite HC-index.
Then for every matrix~$\mW$, the generalized Lyapunov equation~\eqref{eq:genlya} has a solution.
For all solutions~$\mX$ of~\eqref{eq:genlya}, the matrix~$\mE^H \mX\mE$ is unique.
Moreover, if~$\mW$ is positive (semi-)definite, then every solution~$\mX$ of~\eqref{eq:genlya} is positive (semi-)definite on the image of~$P_l$, where~$P_l$ is the spectral projection onto the left deflating subspace associated with the finite eigenvalues of~$\lambda \mE -\mA$.
\end{theorem}
\begin{proof}
For general linear DAE systems with regular matrix pencil $\lambda\mE -\mA$ of DAE-index at most two whose finite eigenvalues  lie in the open left half-plane, the result has been shown in~\cite{Sty02,Sty02a}.
Due to Theorem~\ref{thm:singind}, for semi-dissipative Hamiltonian DAEs~\eqref{DAE:EA} the finite spectrum lies in the closed left half plane, the eigenvalues on the imaginary axis are semi-simple, and the pencil is of DAE-index at most two.
In fact, since the regular matrix pencil has finite HC-index, its finite spectrum lies in the open left half-plane.
Thus, the general results in~\cite{Sty02,Sty02a} directly imply the assertion.
\end{proof}
\begin{remark}
Under the assumptions of Theorem~\ref{thm:genlya}, a solution~$\mX$ of~\eqref{eq:genlya} for a given $\mW\in\PSDHn$ ($\mW\in\PDHn$) yields a weighted semi-norm~$\|\cdot\|_{\tX}$ with $\tX:=\mE^H\mX\mE$ which decays (strict) monotonically along solutions of~\eqref{DAE:EA}, due to Remark~\ref{remark:DAE:HCI}\ref{remark:DAE:HCI:semi-norm:tX}.
To characterize the exponential rate of the dynamic part of the DAE, instead of~\eqref{eq:genlya} one should study the existence of $\mu\in\R$ and $\mX\in\Cnn$ solving a generalized Lyapunov matrix inequality
\begin{equation}\label{DAE:EA:LMI}
\mE^H \mX \mA+ \mA^H \mX\mE
\leq -2\mu \mE^H \mX \mE
\qquad \text{and } \mE^H \mX \mE \in\PSDHn ,
\end{equation}
compare with~\eqref{example:X:LMI} in case of an ODE.
\end{remark}
Using the staircase form~\eqref{staircase:EJR} we obtain a solution procedure for the generalized Lyapunov equation~\eqref{eq:genlya}.

\subsubsection{Solution procedure for a generalized Lyapunov equation}
First, transform the generalized Lyapunov equation~\eqref{eq:genlya} with the unitary matrix~$\mP$ that brings~$\lambda \mE -\mA$ to staircase form~\eqref{staircase:EJR} and via the transformation matrices~$\mL,\mZ$ from Lemma~\ref{lem:reduced} accordingly to
\begin{multline*}
(\mZ^H\mP\mE^H\mP^H\mL^H)(\mL^{-H} \mP \mX\mP^H \mL^{-1}) (\mL \mP \mA \mP^H \mZ)+(\mZ^{H} \mP\mA^H \mP^H \mL^H)(\mL^{-H} \mP \mX\mP^H \mL^{-1}) (\mL\mP \mE \mP^H \mZ)
\\
=-\mZ^H\mP\mE^H \mW \mE\mP^H\mZ \ .
\end{multline*}
Setting~$\mY :=\mL^{-H} \mP \mX\mP^H \mL^{-1}$ and $\widehat \mW :=\mZ^H\mP\mE^H \mW \mE\mP^H\mZ$, the Hermitian matrices~$\mY$ and~$\widehat\mW$ are partitioned analogously to~$\widecheck\mE,\widecheck \mJ,\widecheck \mR$ as $\mY=[\mY_{i,j}]$ and~$\widehat \mW=\widehat \mW^H =[\widehat\mW_{i,j}]$, $i,j=1,\ldots,4$.
Then, the transformed generalized Lyapunov equation reads
\begin{equation}\label{eq:genlya:transformed}
\widehat\mE^H \mY \widehat\mA+ \widehat\mA^H \mY\widehat\mE=-\widehat\mW \,.
\end{equation}
Using the structure of matrices in~\eqref{kcf:almost}, the transformed generalized Lyapunov equation~\eqref{eq:genlya:transformed} reads
\begin{multline*}
\begin{bmatrix}
\widehat\mE_{1,1} & 0 & 0 & 0  \\
0 & \widehat\mE_{2,2} & 0 & 0 \\
0  & 0 & 0 & 0  \\
0  & 0 & 0 & 0
\end{bmatrix}
\begin{bmatrix}
-\mY_{1,4} & \mY_{1,2} \widehat\mA_{2,2} & \mY_{1,3} & \mY_{1,1} \\
-\mY_{2,4} & \mY_{2,2} \widehat\mA_{2,2} & \mY_{2,3} & \mY_{2,1} \\
-\mY_{3,4} & \mY_{3,2} \widehat\mA_{2,2} & \mY_{3,3} & \mY_{3,1} \\
-\mY_{4,4} & \mY_{4,2} \widehat\mA_{2,2} & \mY_{4,3} & \mY_{4,1}
\end{bmatrix}
\\
+
\begin{bmatrix}
-\mY_{1,4} & \mY_{1,2} \widehat\mA_{2,2} & \mY_{1,3} & \mY_{1,1} \\
-\mY_{2,4} & \mY_{2,2} \widehat\mA_{2,2} & \mY_{2,3} & \mY_{2,1} \\
-\mY_{3,4} & \mY_{3,2} \widehat\mA_{2,2} & \mY_{3,3} & \mY_{3,1} \\
-\mY_{4,4} & \mY_{4,2} \widehat\mA_{2,2} & \mY_{4,3} & \mY_{4,1}
\end{bmatrix}^H
\begin{bmatrix}
\widehat\mE_{1,1} & 0 & 0 & 0  \\
0 & \widehat\mE_{2,2} & 0 & 0 \\
0  & 0 & 0 & 0  \\
0  & 0 & 0 & 0
\end{bmatrix}
=-\begin{bmatrix}
\widehat\mW_{1,1} & \widehat\mW_{2,1}^H & 0 & 0  \\
\widehat\mW_{2,1} & \widehat\mW_{2,2} & 0 & 0 \\
0  & 0 & 0 & 0  \\
0  & 0 & 0 & 0
\end{bmatrix} \ ,
\end{multline*}
where we used the identity $\mP\mE^H \mW \mE\mP^H =\mP\mE^H\mP^H\mP \mW \mP^H\mP\mE\mP^H$ to deduce the structure of~$\widehat\mW$.
Multiplying out the left side of~\eqref{eq:genlya:transformed},
reduces the Lyapunov equation to two linear systems
\begin{equation}\label{LS:UL} 
\begin{bmatrix}
\widehat\mE_{1,1} & 0 \\
0 & \widehat\mE_{2,2} \\
\end{bmatrix}
\begin{bmatrix}
-\mY_{1,4} & \mY_{1,2} \widehat\mA_{2,2} \\
-\mY_{2,4} & \mY_{2,2} \widehat\mA_{2,2} \\
\end{bmatrix}
+
\begin{bmatrix}
-\mY_{1,4} & \mY_{1,2} \widehat\mA_{2,2} \\
-\mY_{2,4} & \mY_{2,2} \widehat\mA_{2,2} \\
\end{bmatrix}^H
\begin{bmatrix}
\widehat\mE_{1,1} & 0 \\
0 & \widehat\mE_{2,2} \\
\end{bmatrix}
=-\begin{bmatrix}
\widehat\mW_{1,1} & \widehat\mW_{2,1}^H \\
\widehat\mW_{2,1} & \widehat\mW_{2,2}  \\
\end{bmatrix}
\end{equation}
and
\begin{equation}\label{LS:UR} 
\begin{bmatrix}
\widehat\mE_{1,1} &0  \\
0 &\widehat\mE_{2,2}
\end{bmatrix}
\begin{bmatrix}
\mY_{1,3} &\mY_{1,1}  \\
\mY_{2,3} &\mY_{2,1}
\end{bmatrix}
=0 \ ,
\end{equation}
corresponding to the upper left and upper right block, respectively.
Since the matrix~$\diag(\widehat\mE_{1,1},\widehat\mE_{2,2})$ is invertible, equation~\eqref{LS:UR} implies that the blocks $\mY_{1,3}=\mY_{3,1}^H,\mY_{1,1},\mY_{2,3}=\mY_{3,2}^H,\mY_{2,1}=\mY_{1,2}^H$ are~$0$.
The blocks $\mY_{3,3},\mY_{4,3}=\mY_{3,4}^H,\mY_{4,4}$ do not occur in the equations and can therefore be chosen arbitrarily.
Using~$\mY_{2,1}=\mY_{1,2}^H=0$, equation~\eqref{LS:UL} simplifies to
\begin{subequations}
\begin{align}
\widehat\mE_{2,2} \mY_{2,2} \widehat\mA_{2,2} +\widehat\mA_{2,2}^H \mY_{2,2} \widehat\mE_{2,2}
&=-\widehat \mW_{2,2} \ , \label{eq:dynlyaeq} \\
-\widehat\mE_{1,1}\mY_{1,4} -\mY_{1,4}^H\widehat\mE_{1,1}
&=-\widehat \mW_{1,1} \ , \label{eq:y14} \\
-\mY_{2,4}^H \widehat\mE_{2,2}
&=-\widehat\mW_{2,1}^H \label{eq:y24:implicit} \ ,
\end{align}
\end{subequations}
which are independent equations for $\mY_{2,2}$, $\mY_{1,4}$ and $\mY_{2,4}$, respectively.
Equation~\eqref{eq:dynlyaeq} is associated with the dynamic part of the system.
Since~$\widehat\mE_{2,2}$ is positive definite and~$\widehat\mA_{2,2}$ is negative hypocoercive (i.e. all eigenvalues have negative real part) under the assumptions of Theorem~\ref{thm:genlya}, the Lyapunov equation~\eqref{eq:dynlyaeq} has a unique positive (semi-)definite solution~$\mY_{2,2}$ for every positive (semi-)definite~$\widehat\mW_{2,2}$, see e.g.~\cite{Sty02,Sty02a}.
Next, since $\widehat\mE_{1,1}$ is Hermitian positive definite, the Lyapunov equation~\eqref{eq:y14} has a unique positive (semi-)definite solution~$\mY_{1,4}$ for every positive (semi-)definite~$\widehat\mW_{1,1}$.
Finally, using~\eqref{eq:y24:implicit} and $\widehat\mE_{2,2}>0$ yields
\begin{equation}\label{eq:y24}
\mY_{2,4}^H
=\widehat\mW_{2,1}^H \widehat\mE_{2,2}^{-1} \ .
\end{equation}
Altogether, the solution set of~\eqref{eq:genlya} consists of all Hermitian matrices of the form
\[
\mX
= \mP^H \mL^H
\begin{bmatrix}
0 & 0 & 0 & \mY_{1,4} \\
0 & \mY_{2,2} & 0 & \mY_{2,4} \\
0 & 0 & \mY_{3,3} & \mY_{3,4} \\
\mY_{1,4} & \mY_{2,4}^H & \mY_{3,4}^H & \mY_{4,4}
\end{bmatrix} \mL \mP,
\]
with arbitrary blocks $\mY_{3,3},\mY_{4,3}=\mY_{3,4}^H,\mY_{4,4}$, positive (semi-)definite Hermitian matrices~$\mY_{2,2}$ and~$\mY_{1,4}$ satisfying~\eqref{eq:dynlyaeq} and~\eqref{eq:y14}, respectively, and finally~$\mY_{2,4}$ determined by~\eqref{eq:y24}.

\section{Examples}
\label{sec:example}

In this section we  illustrate the above exposition on some simple ODE and PDE examples.

\subsection{Stokes equation on the 2D torus} \label{ssec:Stokes}
First we consider the time-dependent, incompressible Stokes equation of fluid dynamics on the 2D torus $\mathbb T^2 := (0,2\pi)^2$,
\begin{equation} \label{eq:Stokes}
\begin{cases}
u_t +\nabla p &=\nu \Delta u\,, \quad t>0\,, \\
\diver u &=0\,,
\end{cases}
\end{equation}
for the vector-valued velocity field~$u=u(x,t)$ and the scalar pressure~$p =p(x,t)$ in the space variable~$x\in\T^2$ and the time variable~$t\geq 0$.
The constant $\nu>0$ denotes the viscosity coefficient.
Due to the periodic boundary conditions in~\eqref{eq:Stokes}, this model actually could be simplified right away:
Taking the divergence of the first equation in~\eqref{eq:Stokes} yields
$\Delta p(\cdot,t) =0$ and hence~$p(\cdot,t)$ is constant in~$x$.
It also shows that the vector-valued heat equation $u_t =\nu \Delta u$ preserves the incompressibility if the initial condition satisfies $\diver u(0)=0$, which is assumed in the sequel.
But to illustrate negative hypocoercivity of matrix pencils in semi-dissipative Hamiltonian DAEs we shall ignore this possible simplification and rather follow our discussion from~\S\ref{sec:DAE}.

Due to the periodic setting we consider the Fourier expansion of~\eqref{eq:Stokes} with
\[
u(x,t)
= \sum_{k\in\Z^2} u_k(t) e^{i k\cdot x} \,, \qquad
p(x,t)
= \sum_{k\in\Z^2} p_k(t) e^{i k\cdot x} \,.
\]
The Fourier coefficients~$u_k(t)\in\C^2$, $p_k(t)\in\C$, $k\in\Z^2$, satisfy the decoupled evolution equations
\begin{equation}\label{eq:Stokes:Fourier}
\begin{cases}
\ddt u_k &= -i k p_k -\nu |k|^2 u_k \ , \quad t>0 \ , \\
i k\cdot u_k &= 0 \ .
\end{cases}
\end{equation}
The mode~$k=0$ satisfies $u_0(t)=$ const. (corresponding to momentum conservation) and~$p_0(t)=$ arbitrary.
To enforce unique solvability of~\eqref{eq:Stokes}, we normalize the pressure as $p_0(t)\equiv 0$.
For~$k\ne 0$ we write~\eqref{eq:Stokes:Fourier} as a system of decoupled DAEs, each having DAE-index~$2$:
\begin{equation} \label{DAE:Stokes}
 \mE \dot{w}_k (t) =\mA_k \xk \ , \quad t\geq 0\ ,
\end{equation}
for~$\xk :=[u_k^1,u_k^2,p_k]^\top \in\C^3$ with the matrices~$\mE:=\diag(1,1,0)$ and
\begin{equation} \label{DAE:Stokes:A}
\mA_k :=
 \begin{bmatrix}
  -\nu |k|^2 & 0 & -i k_1 \\
  0 & -\nu |k|^2 & -i k_2 \\
  -i k_1 & -i k_2 & 0
 \end{bmatrix} \ .
\end{equation}
The modal functions~$\xk(t)$, $k\in\Z^2$ correspond to the function~$x(t)$ in~\S\ref{sec:intro}--\ref{sec:DAE}, since~$x=[x_1, x_2]^\top$ is used here for the spatial variable.
Following the notation from~\S\ref{sec:intro}, we  decompose $\mA_k$ as $\mA_k =\mJ_k -\mR_k$ with $\mR_k :=\diag(\nu |k|^2,\nu |k|^2, 0)$ and
\begin{equation} \label{DAE:Stokes:J}
\mJ_k :=
 \begin{bmatrix}
  0 & 0 & -i k_1 \\
  0 & 0 & -i k_2 \\
  -i k_1 & -i k_2 & 0
 \end{bmatrix} \ .
\end{equation}
In order to define the HC-index of~\eqref{DAE:Stokes} we  transform~\eqref{DAE:Stokes} to staircase form:
A straightforward application of the Staircase Algorithm in Lemma~\ref{lem:tSF} to the triple~$(\mE,\mJ_k,\mR_k)$ yields the (unitary) congruence transformation
\[
\widecheck \mJ_k
=\mP_k \mJ_k \mP_k^H
=\begin{bmatrix}
  0 & 0 & -|k| \\
  0 & 0 & 0 \\
  |k| & 0 & 0
 \end{bmatrix} \ , \quad
\widecheck \mR_k
=\mP_k \mR_k \mP_k^H
=\mR_k \ , \quad
\widecheck \mE
=\mP_k \mE \mP_k^H
=\mE \ , \quad
\]
with constants $n_1 =n_2 =n_4=1$, $n_3 =n_5 =0$ and
\begin{equation} \label{ex:Stokes:Pk}
\mP_k
=\tfrac1{|k|}
 \begin{bmatrix}
  k_1 & k_2 & 0 \\
  -k_2 & k_1 & 0 \\
  0 & 0 & i|k|
 \end{bmatrix} \ .
\end{equation}
The evolution of~\eqref{DAE:Stokes} translates via $y_k =[y_{k,1},y_{k,2},y_{k,4}]^\top :=\mP_k \xk$ (note that $y_{k,3}$ is void) into the staircase form
\begin{equation} \label{DAE:Stokes:SF}
 \widecheck \mE \dot y_k (t) =(\widecheck \mJ_k -\widecheck \mR_k) y_k(t) \ , \quad t\geq 0 \ .
\end{equation}

\smallskip
From~\ref{P:n5}, \ref{P:DAE-index} in Corollary~\ref{cor:index} we see that all pencils $\lambda \widecheck \mE -(\widecheck \mJ_k -\widecheck \mR_k)$, $k\ne 0$ are regular of DAE-index~$2$.
Similarly to~\eqref{DAE:EA:UC:coordinates} we obtain
\[
 y_{k,1} =\tfrac1{|k|} k\cdot u_k =0 \,, \qquad
 y_{k,4} = i p_k =0 \qquad
 \text{for all } k\ne 0 .
\]
Following~\ref{EJR:ODE}, the evolution can be reduced to the dynamic part:
\[
\mE_{2,2} \dot y_{k,2}
=
\widehat \mA_{k,2,2}  y_{k,2} \ ,
\quad
t\geq 0 \ ,
\]
with $\mE_{2,2}=1$, $\widehat \mA_{k,2,2} =-\nu |k|^2$.
Hence, the evolution reduces to
\begin{equation} \label{DAE:Stokes:y2}
 \dot y_{k,2} =-\nu |k|^2 y_{k,2} \ , \quad k\ne 0
\end{equation}
with HC-index $0$, as introduced in Definition~\ref{def:DAE:mHC}.
Equation~\eqref{DAE:Stokes:y2} is the modal decomposition of the (dissipative) heat equation on $\mathbb T^2$.
Hence, the solution $(u(\cdot,t),p(\cdot,t))$ of the Stokes equation~\eqref{eq:Stokes} converges, as $t\to \infty$, to the constant equilibrium $(u_0,p_0)$ with the exponential decay rate~$\nu=\min_{k\ne 0} (\nu|k|^2)$.

\smallskip
We remark that the same analysis carries over to the time-dependent Oseen equation~\cite{BaJi08}
\begin{equation} \label{eq:Oseen}
\begin{cases}
u_t +(b\cdot\nabla)u +\nabla p &=\nu \Delta u\,, \quad t>0\,, \\
\diver u &=0\,,
\end{cases}
\end{equation}
on~$\mathbb T^d$ with some constant~$b\in\R^d$.

When modifying~\eqref{eq:Oseen} with $d=2$ into an anisotropic Oseen equation with viscosity only in the $x_2$--direction, the dynamics becomes more interesting:
For a constant convection field~$b$ in the $x_1$--direction, the modes still decouple but the generator of the evolution is neither coercive nor hypocoercive.
If the convection field is non-constant, e.g. $b=[\sin x_2,\,0]^\top$, the spatial modes are coupled and the generator of the (infinite-dimensional) problem becomes hypocoercive.
A detailed analysis will be the topic of a forthcoming paper.

\subsection{Network of gas pipelines}
\label{ssec:network}

We consider a simple model for acoustic waves in a fluid flow through a network of gas pipelines, which has been studied in~\cite{BMXZ18,EK18,EKLSMM18}:
Since the cross-section of a single pipe is usually much smaller than the length of the pipe, single pipes are modeled as one-dimensional.
Then, a network of pipelines is represented as a finite directed and connected graph~$\cG(\cV,\cE)$ with vertices $v\in\cV$ and edges $e\in\cE$.

\textit{Differential equations.}
On every edge~$e$ of the graph ($\cong$ a pipe of the physical network), the fluid flow is modeled via its pressure $p^e =p^e(x,t)$ and mass flux~$q^e =q^e(x,t)$.
The acoustic pressure wave is then subject to the linear damped hyperbolic system
\begin{subequations}\label{Gas:DE}
 \begin{align}
  a^e \partial_t p^e +\partial_x q^e &= 0
  &&\text{on } e\in\cE,\ x\in[0,\ell^e],\ t>0, \\
  b^e \partial_t q^e +\partial_x p^e &= -d^e q^e
  &&\text{on } e\in\cE,\ x\in[0,\ell^e],\ t>0,
 \end{align}
\end{subequations}
where the parameters $a^e, b^e, d^e, \ell^e$ encode properties of fluid and pipe, and are assumed to be positive and constant on each pipe/edge.
In particular, the parameter~$d^e>0$ is related to damping due to the friction at the pipe walls.
For a derivation of a nonlinear variant of~\eqref{Gas:DE} from Euler equations see, e.g., \S3.2.2 in~\cite{BrGaHe11}.
\new{More generally, for the reformulation of hyperbolic conservation laws in symmetric form we refer to~\cite{FL71}.}

\textit{Algebraic constraints.}
The fluid flow at a junction is assumed to (i) conserve mass and (ii) exhibit a unique pressure, which translates into the frequently used coupling conditions:
Consider the subset of all inner vertices $\cV_0\subset\cV$,
\begin{subequations}\label{Gas:coupling}
 \begin{align}
  \sum_{e\in\cE(v)} n^e(v) q^e(v) &= 0
  &&\text{for } v\in\cV_0,\ t>0, \label{coupling:mass} \\
  p^e(v) &= p^{e'}(v)
  &&\text{for } v\in\cV_0,\ e,e'\in\cE(v),\ t>0,
 \end{align}
\end{subequations}
where $n^e(v)=\pm1$ depending whether pipe~$e$ starts or ends at~$v$, $\cE(v)$ are the edges adjacent to vertex~$v$, and $q^e(v), p^e(v)$ denote the respective functions evaluated at the vertex~$v$ but still depending on time.

\textit{Boundary conditions.}
At the boundary vertices $v\in\cV_\partial :=\cV\setminus\cV_0$ ($\cong$ ports of the network), we set
\begin{align}\label{Gas:input}
 p^e(v) &=0
 &&\text{for } v\in\cV_\partial,\ e\in\cE(v),\ t>0,
\end{align}
i.e. homogeneous boundary conditions.
In~\cite{BMXZ18,EK18,EKLSMM18}, the system is controlled via the (given) pressure~$p^e(v) =u_v$ at the port $v\in\cV_\partial$.

\textit{Initial conditions.}
The specification of the model is completed by assuming knowledge of the initial conditions
\begin{align}\label{Gas:IC}
 p^e(0) =p_0, \ \qquad
 q^e(0) &=q_0, \ \qquad
 &&\text{on } \cE .
\end{align}

The partial differential-algebraic system~\eqref{Gas:DE}--\eqref{Gas:input} encodes several interesting properties which are directly related to the underlying physical principles:
\begin{enumerate}[(N1)]
 \item \label{N1}
 \emph{Global conservation of mass.}
 \item \label{N2}
 \emph{A port-Hamiltonian structure.}
 \item \label{N3}
 \emph{Exponential stability and convergence to equilibrium.}
 More precisely, 
 the energy of the system
\[
E(t)
:=\tfrac12 \sum_{e\in\cE} \int_e a^e|p^e|^2 +b^e |q^e|^2 \d[x]\; \text{  decays exponentially  }\;
E(t)
\leq C e^{-\gamma t} E(0),
\quad %
t\geq 0,
\]
 with constants $C$ and $\gamma$ that are independent of the particular solution.
 %
 \item \label{N4}
 \emph{Unique steady state~$(\bar p,\bar q)=(0,0)$} for the corresponding stationary problem.
\end{enumerate}

We follow here the presentation in Example~11 and Example~24 of~\cite{BMXZ18} without input control and without output:
A mixed finite element discretization that preserves the structural properties~\ref{N1}--\ref{N4} leads to a block structured constant coefficient port-Hamiltonian DAE system
\begin{equation} \label{pHDAE:network:0}
\mE \dot x
= (\mJ -\mR) x\ ,
\qquad %
x(0)
= x^0 \ ,
\end{equation}
with
\begin{align*}
\mE
&=\begin{bmatrix}
  \mM_1 & 0 & 0 \\
  0 & \mM_2 & 0 \\
  0 & 0 & 0
 \end{bmatrix} ,
& %
\mJ
&=\begin{bmatrix}
  0 & -\tG & 0 \\
  \tG^\top & 0 & \tN^\top \\
  0 & -\tN & 0
 \end{bmatrix} ,
& %
\mR
&=\begin{bmatrix}
  0 & 0 & 0 \\
  0 & \tD & 0 \\
  0 & 0 & 0
 \end{bmatrix} ,
& %
x
&:=\begin{bmatrix}
  x_1 \\ x_2 \\ x_3
 \end{bmatrix} ,
\end{align*}
where the vector valued functions $x_1:\R\to\R^{m_1}$, $x_2:\R\to\R^{m_2}$, represent the discretized pressure and flux, respectively, and $x_3:\R\to\R^{m_3}$, represents the Lagrange multiplier for satisfying the space-discretized constraints~\eqref{coupling:mass}.
Here we redefine the variable~$x$ to relate to our standard notation in~\S\ref{sec:intro}--\ref{sec:DAE}.
The coefficient matrices $\mM_1\in\R^{m_1\times m_1}$, $\mM_2 \in\R^{m_2 \times m_2}$ and $\tD  \in\R^{m_2 \times m_2}$ are symmetric and positive definite.
Moreover $m_1 +m_3 \leq m_2$, and the matrix $[\tG^\top,\ \tN^\top]\in\R^{m_2\times(m_1 +m_3)}$ has trivial null-space.
The Hamiltonian is given by $\mathsf H(x) :=\tfrac12 (\mE x)^\top x =\tfrac12 (x_1^\top \mM_1 x_1 +x_2^\top \mM_2 x_2)$.
The system has DAE-index~$2$, see~\cite{EKLSMM18,KunM06}.

\medskip
In the following, we study the negative hypocoercivity of the matrix pencil associated to the semi-dissipative Hamiltonian system~$\mE \dot x =(\mJ -\mR)x$.
To derive the staircase form for the triple~$(\mE,\mJ,\mR)$,
we follow the proof of Lemma~\ref{lem:tSF}:
The matrices $\mE,\mJ,\mR\in\Rnn$, $n=m_1+m_2+m_3$ are already in the desired form of Step~1 and Step~2 such that $\mP=\mI$ with $\tn_1 =m_1 +m_2$ and $\tn_2 =0$.
In Step~3, an SVD of $\tJ_{3,1} =[0, {-\tN}]\in\R^{\tn_3\times \tn_1}$, $\tn_3 =m_3$ is performed.
Since~$\tN\in\R^{m_3\times m_2}$ with $m_3\leq m_2$ has full row rank, the matrix~$\tN$ has an SVD of the form
\begin{equation}\label{SVD:tN}
\tN =\mU \begin{bmatrix} \Sigma & 0 \end{bmatrix} \mV^\top ,
\end{equation}
with real orthogonal matrices $\mU$, $\mV$ and a regular diagonal matrix~$\Sigma\in\R^{m_3\times m_3}$,
see also~\cite{EKLSMM18}.
Thus, an SVD of~$\tJ_{3,1}$ reads
\[
 \tJ_{3,1}
 =
 \mU_{3,1}
 \left[\begin{array}{cc|c} \Sigma & 0 & 0 \end{array}\right]
 \mV_{3,1}^\top , \qquad
 \text{where }
 \mU_{3,1}
 =-\mU , \quad
 \mV_{3,1}
 =
 \begin{bmatrix} 0 & \mI_{m_1} \\ \mV & 0 \end{bmatrix} .
\]
Following Lemma~\ref{lem:tSF}, the real orthogonal matrix
$\mP =\mP_3 =\diag(\mV_{3,1}^\top , \mU_{3,1}^\top)$ yields
\begin{align*}
\widecheck\mE
&=\begin{bmatrix}
  \mV^\top \mM_2 \mV & 0 & 0 \\
  0 & \mM_1 & 0 \\
  0 & 0 & 0
 \end{bmatrix} ,
& %
\widecheck\mJ
&=\begin{bmatrix}
  0 & \mV^\top \tG^\top & -\mV^\top\tN^\top\mU \\
  -\tG \mV & 0 & 0 \\
  \mU^\top\tN\mV & 0 & 0
 \end{bmatrix} ,
& %
\widecheck\mR
&=\begin{bmatrix}
  \mV^\top \tD \mV & 0 & 0 \\
  0 & 0 & 0 \\
  0 & 0 & 0
 \end{bmatrix} .
\end{align*}
Using~\eqref{SVD:tN} we define
\begin{align}\label{pHDAE:network:EJR}
\widecheck\mE
&=:
\left[\begin{array}{cc|c|c}
\mM_{1,1} & \mM_{2,1}^\top & 0 & 0 \\
\mM_{1,2} & \mM_{2,2} & 0 & 0 \\
\hline
0 & 0 & \mM_1 & 0 \\
\hline
0 & 0 & 0 & 0
\end{array}\right] ,
& %
\widecheck\mJ
&=:
\left[\begin{array}{cc|c|c}
0 & 0 & \mG_{1,1}^\top & -\Sigma \\
0 & 0 & \mG_{1,2}^\top & 0 \\
\hline
-\mG_{1,1} & -\mG_{1,2} & 0 & 0 \\
\hline
\Sigma & 0 & 0 & 0
\end{array}\right] ,
& %
\widecheck\mR
&=:
\left[\begin{array}{cc|c|c}
\mD_{1,1} & \mD_{2,1}^\top & 0 & 0 \\
\mD_{2,1} & \mD_{2,2} & 0 & 0 \\
\hline
0 & 0 & 0 & 0 \\
\hline
0 & 0 & 0 & 0
\end{array}\right] \ ,
\end{align}
where the lines indicate the previous partitioning (but not the partitioning of~\eqref{staircase:EJR}).
For the partitioning of~\eqref{staircase:EJR}, we find $n_1=n_4=m_3$, $n_2=n-2m_3$, and $n_3=n_5=0$.
Following Corollary~\ref{cor:index} (or by direct reasoning), we identify the underlying implicit ODE~\eqref{underlyingODE:y} again as
$\mE_{2,2} \dot y_2 =(\widehat \mJ_{2,2}-\widehat \mR_{2,2}) y_2$ with $y_2\in\R^{n_2}$ and
\begin{align}\label{GAS:underlyingODE}
\mE_{2,2}
&=\begin{bmatrix}
  \mM_{2,2} & 0 \\
  0 & \mM_1
 \end{bmatrix} ,
& %
\widehat\mJ_{2,2}
&=\begin{bmatrix}
  0 & \mG_{1,2}^\top \\
  -\mG_{1,2} & 0
 \end{bmatrix} ,
& %
\widehat\mR_{2,2}
&=\begin{bmatrix}
  \mD_{2,2} & 0 \\
  0 & 0
 \end{bmatrix} .
\end{align}
To determine the HC-index of this implicit ODE (and of the original DAE) we have to use a congruence transformation to an ODE as in~\eqref{ODE:A-tilde}, where $\tA
 =(\mE_{2,2}^{1/2})^{-1} (\widehat \mJ_{2,2}-\widehat \mR_{2,2}) (\mE_{2,2}^{1/2})^{-1} =:\tJ-\tR$ with
\begin{align}\label{Gas:ODE:tJ-tR}
\tJ
&=\begin{bmatrix}
  0 & (\mM_{2,2}^{1/2})^{-1} \mG_{1,2}^\top (\mM_1^{1/2})^{-1} \\
  -(\mM_1^{1/2})^{-1} \mG_{1,2} (\mM_{2,2}^{1/2})^{-1} & 0
 \end{bmatrix} ,
& %
\tR
&=\begin{bmatrix}
   (\mM_{2,2}^{1/2})^{-1} \mD_{2,2} (\mM_{2,2}^{1/2})^{-1} & 0 \\
  0 & 0
 \end{bmatrix} .
\end{align}
The matrix~$\tR$ is Hermitian positive semi-definite, but
\[
\tR +\tJ\ \tR\ \tJ^\top
=
\begin{bmatrix}
 (\mM_{2,2}^{1/2})^{-1} \mD_{2,2} (\mM_{2,2}^{1/2})^{-1} & 0 \\
 0 & -(\mM_1^{1/2})^{-1} \mG_{1,2} \mM_{2,2}^{-1} \mD_{2,2} \mM_{2,2}^{-1} \mG_{1,2}^\top (\mM_1^{1/2})^{-1}
\end{bmatrix}
\]
is Hermitian positive definite, since $\mD_{2,2}, \mM_{2,2}, \mM_1$ are Hermitian positive definite (as principal minors of Hermitian positive definite matrices) and $\mG_{1,2}^\top \in\R^{m_1 \times (m_2-m_3)}$ has trivial null-space.

Thus, the ODE~\eqref{ODE:A-tilde} with $\tA=\tJ-\tR$ has HC-index~1.
Due to Corollary~\ref{cor:hypodae}, the origin~$y_2=0$ is asymptotically stable for ODE~\eqref{underlyingODE:y} with~\eqref{GAS:underlyingODE}.
Hence, for every consistent initial condition the solution of DAE~\eqref{DAE:EA:UnitarilyCongruent} with~\eqref{pHDAE:network:EJR} converges exponentially to the unique steady state~0.
Thus, we verify again that this specific mixed finite element discretization preserves the properties~\ref{N3} and~\ref{N4}.

\section{Conclusions}\label{sec:conclusions}
We have studied linear ODEs and DAEs exhibiting hypocoercivity and related this concept to classical concepts from control theory.
\new{In particular we showed for ODEs that the hypocoercivity index can be obtained in a numerically stable manner from the staircase form of the ODE-generator matrix. 
For DAEs we extended the notion of hypocoercivity index to its matrix pencil, and we proved that this index still characterizes the short time behavior of DAE-solutions pertaining to consistent initial conditions.} 
The results are illustrated via two infinite-dimensional application problems.

\section*{Acknowledgments}

The first author (FA) was supported by the FWF-funded SFB \# F65.
The second author (AA) was partially supported by the FWF-doctoral school ``Dissipation and dispersion in non-linear partial differential equations'' and the FWF-funded SFB \# F65.
The third author (VM) was supported by DFG SFB \# 910.

\bibliographystyle{plain}

\fontfamily{\rmdefault}
\fontsize{8bp}{10}
\selectfont\baselineskip=10pt
\bibliography{AAM-2020-HCIndex}

\begin{thebibliography}{10}

\bibitem{AAC16}
F.~Achleitner, A.~Arnold, and E.~A. Carlen.
\newblock On linear hypocoercive {BGK} models.
\newblock In {\em From particle systems to partial differential equations.
  {III}}, volume 162 of {\em Springer Proc. Math. Stat.}, pages 1--37.
  Springer, Cham, 2016.

\bibitem{AAC18}
F.~Achleitner, A.~Arnold, and E.~A. Carlen.
\newblock On multi-dimensional hypocoercive {BGK} models.
\newblock {\em Kinet. Relat. Models}, 11(4):953--1009, 2018.

\bibitem{AAC20}
F.~Achleitner, A.~Arnold, and E.~A. Carlen.
\newblock The hypocoercivity index for the short- and large-time behavior of
  {ODE}s.
\newblock {\em work in progress}, 2021.

\bibitem{AAM21}
F.~Achleitner, A.~Arnold, and V.~Mehrmann.
\newblock The hypocoercivity index in the infinite-dimensional setting.
\newblock {\em work in progress}, 2021.

\bibitem{AAS19}
F.~Achleitner, A.~Arnold, and B.~Signorello.
\newblock On optimal decay estimates for {ODE}s and {PDE}s with modal
  decomposition.
\newblock In {\em Stochastic dynamics out of equilibrium}, volume 282 of {\em
  Springer Proc. Math. Stat.}, pages 241--264. Springer, Cham, 2019.

\bibitem{AASt15}
F.~Achleitner, A.~Arnold, and D.~St\"{u}rzer.
\newblock Large-time behavior in non-symmetric {F}okker--{P}lanck equations.
\newblock {\em Riv. Math. Univ. Parma (N.S.)}, 6(1):1--68, 2015.

\bibitem{Adr95}
L.~Ya. Adrianova.
\newblock {\em Introduction to linear systems of differential equations}.
\newblock Trans. Math. Monographs, Vol. 146, AMS, Providence, RI, 1995.

\bibitem{ArEr14}
A.~Arnold and J.~Erb.
\newblock Sharp entropy decay for hypocoercive and non-symmetric
  {F}okker--{P}lanck equations with linear drift.
\newblock {\em arXiv preprint arXiv:1409.5425}, 2014.

\bibitem{ArJiWo19}
A.~Arnold, S.~Jin, and T.~W\"ohrer.
\newblock Sharp decay estimates in local sensitivity analysis for evolution
  equations with uncertainties: from odes to linear kinetic equations.
\newblock {\em arXiv preprint arXiv:1904.01190}, 2019.

\bibitem{ASS20}
A.~Arnold, C.~Schmeiser, and B.~Signorello.
\newblock Propagator norm and sharp decay estimates for {F}okker--{P}lanck
  equations with linear drift.
\newblock {\em arXiv preprint arXiv:2003.01405}, 2020.

\bibitem{AS21}
A.~Arnold and B.~Signorello.
\newblock Optimal non-symmetric {Fokker--Planck} equation for the convergence
  to a given equilibrium.
\newblock {\em arXiv preprint arXiv:2106.15742}, 2021.

\bibitem{BaJi08}
H.-O. Bae and B.~J. Jin.
\newblock Estimates of the wake for the 3{D} {O}seen equations.
\newblock {\em Discrete Contin. Dyn. Syst. Ser. B}, 10(1):1--18, 2008.

\bibitem{BaVu90}
C.~J.~K. Batty and Q.~P. V\~{u}.
\newblock Stability of individual elements under one-parameter semigroups.
\newblock {\em Trans. Amer. Math. Soc.}, 322(2):805--818, 1990.

\bibitem{BeMeVD19}
C.~Beattie, V.~Mehrmann, and P.~Van~Dooren.
\newblock Robust port-{H}amiltonian representations of passive systems.
\newblock {\em Automatica J. IFAC}, 100:182--186, 2019.

\bibitem{BMXZ18}
C.~Beattie, V.~Mehrmann, H.~Xu, and H.~Zwart.
\newblock Linear port-{H}amiltonian descriptor systems.
\newblock {\em Math. Control Signals Systems}, 30(4):Art. 17, 27, 2018.

\bibitem{Be18}
D.~S. Bernstein.
\newblock {\em Scalar, vector, and matrix mathematics}.
\newblock Princeton University Press, Princeton, NJ, 2018.

\bibitem{scalapack}
L.~S. Blackford, J.~Choi, A.~Cleary, E.~D'Azevedo, J.~Demmel, I.~Dhillon,
  J.~Dongarra, S.~Hammarling, G.~Henry, A.~Petitet, K.~Stanley, D.~Walker, and
  R.~C. Whaley, editors.
\newblock {\em {ScaLAPACK} Users' Guide}.
\newblock Software, Environments and Tools. Society for Industrial and Applied
  Mathematics, Philadelphia, PA, USA, 1997.

\bibitem{BrGaHe11}
J.~Brouwer, I.~Gasser, and M.~Herty.
\newblock Gas pipeline models revisited: model hierarchies, nonisothermal
  models, and simulations of networks.
\newblock {\em Multiscale Model. Simul.}, 9(2):601--623, 2011.

\bibitem{BruM07}
T.~{Br\"ull} and V.~{Mehrmann}.
\newblock {STCSSP}: {A} {FORTRAN~77} routine to compute a structured staircase
  form for a (skew-)symmetric/{(skew-)}symmetric pencil.
\newblock Preprint 31-2007, Institut f\"ur Mathematik, TU Berlin, 2007.

\bibitem{ByeGM97}
R.~{Byers}, T.~{Geerts}, and V.~{Mehrmann}.
\newblock Descriptor systems without controllability at infinity.
\newblock {\em SIAM J. Control Optim.}, 35(2):462--479, 1997.

\bibitem{ByeMX07}
R.~{Byers}, V.~{Mehrmann}, and H.~{Xu}.
\newblock A structured staircase algorithm for skew-symmetric/symmetric
  pencils.
\newblock {\em Electron. Trans. Numer. Anal.}, 26:1--33, 2007.

\bibitem{Da04}
B.~N. Datta.
\newblock {\em Numerical methods for linear control systems}.
\newblock Elsevier Academic Press, San Diego, CA, 2004.

\bibitem{DieRV97}
L.~Dieci, R.~D. Russell, and E.~S. Van~Vleck.
\newblock On the computation of {L}yapunov exponents for continuous dynamical
  systems.
\newblock {\em SIAM J. Numer. Anal.}, 34(1):402--423, 1997.

\bibitem{DieV02b}
L.~Dieci and E.~S. {Van Vleck}.
\newblock Lyapunov and other spectra: a survey.
\newblock In {\em Collected lectures on the preservation of stability under
  discretization (Fort Collins, CO, 2001)}, pages 197--218. SIAM, Philadelphia,
  PA, 2002.

\bibitem{DoMoSc09}
J.~Dolbeault, C.~Mouhot, and C.~Schmeiser.
\newblock Hypocoercivity for kinetic equations with linear relaxation terms.
\newblock {\em C. R. Math. Acad. Sci. Paris}, 347(9-10):511--516, 2009.

\bibitem{DoMoSc15}
J.~Dolbeault, C.~Mouhot, and C.~Schmeiser.
\newblock Hypocoercivity for linear kinetic equations conserving mass.
\newblock {\em Trans. Amer. Math. Soc.}, 367(6):3807--3828, 2015.

\bibitem{DuLM13}
N.~H. Du, V.~H. Linh, and V.~Mehrmann.
\newblock Robust stability of differential-algebraic equations.
\newblock In A.~Ilchmann and T.~Reis, editors, {\em Surveys in
  Differential-Algebraic Equations I}, pages 63--95. Springer, Berlin,
  Heidelberg, 2013.

\bibitem{EK18}
H.~Egger and T.~Kugler.
\newblock Damped wave systems on networks: exponential stability and uniform
  approximations.
\newblock {\em Numer. Math.}, 138(4):839--867, 2018.

\bibitem{EKLSMM18}
H.~Egger, T.~Kugler, B.~Liljegren-Sailer, N.~Marheineke, and V.~Mehrmann.
\newblock On structure-preserving model reduction for damped wave propagation
  in transport networks.
\newblock {\em SIAM J. Sci. Comput.}, 40(1):A331--A365, 2018.

\bibitem{EmmM13}
E.~{Emmrich} and V.~{Mehrmann}.
\newblock Operator differential-algebraic equations arising in fluid dynamics.
\newblock {\em Comput. Methods Appl. Math.}, 13(4):443--470, 2013.

\bibitem{FL71}
K.~O. Friedrichs and P.~D. Lax.
\newblock Systems of conservation equations with a convex extension.
\newblock {\em Proc. Nat. Acad. Sci. U.S.A.}, 68:1686--1688, 1971.

\bibitem{GaMi13}
S.~Gadat and L.~Miclo.
\newblock Spectral decompositions and {$L^2$}-operator norms of toy
  hypocoercive semi-groups.
\newblock {\em Kinet. Relat. Models}, 6(2):317--372, 2013.

\bibitem{Gan59a}
F.~R. Gantmacher.
\newblock {\em The theory of matrices. {V}ols. 1, 2}.
\newblock Translated by K. A. Hirsch. Chelsea Publishing Co., New York, 1959.

\bibitem{GraMQSW16}
N.~Gr{\"a}bner, V.~Mehrmann, S.~Quraishi, C.~Schr\"oder, and U.~{von W}agner.
\newblock Numerical methods for parametric model reduction in the simulation of
  disc brake squeal.
\newblock 96(DOI: 10.1002/zamm.201500217):1388--1405, 2016.

\bibitem{GuMo16}
A.~Guillin and P.~Monmarch\'{e}.
\newblock Optimal linear drift for the speed of convergence of an hypoelliptic
  diffusion.
\newblock {\em Electron. Commun. Probab.}, 21:Paper No. 74, 14, 2016.

\bibitem{GuMo16E}
A.~Guillin and P.~Monmarch\'{e}.
\newblock Erratum: {O}ptimal linear drift for the speed of convergence of an
  hypoelliptic diffusion.
\newblock {\em Electron. Commun. Probab.}, 22:Paper No. 15, 2, 2017.

\bibitem{HiPr10}
D.~Hinrichsen and A.~J. Pritchard.
\newblock {\em Mathematical systems theory {I}}.
\newblock Springer, Heidelberg, 2010.

\bibitem{HoJo13}
R.~A. Horn and C.~R. Johnson.
\newblock {\em Matrix analysis}.
\newblock Cambridge University Press, Cambridge, second edition, 2013.

\bibitem{JoSm05}
C.~R. Johnson and R.~L. Smith.
\newblock Closure properties.
\newblock In F.~Zhang, editor, {\em The Schur Complement and Its Applications},
  pages 111--136. Springer US, Boston, MA, 2005.

\bibitem{JoSm06}
C.~R. Johnson and R.~L. Smith.
\newblock Closure of matrix classes under {S}chur complementation, including
  singularities.
\newblock In {\em Algebra and its applications}, volume 419 of {\em Contemp.
  Math.}, pages 185--200. Amer. Math. Soc., Providence, RI, 2006.

\bibitem{KunM06}
P.~Kunkel and V.~Mehrmann.
\newblock {\em Differential-algebraic equations}.
\newblock European Mathematical Society (EMS), Z\"{u}rich, 2006.

\bibitem{Le11}
W.~S. Levine, editor.
\newblock {\em The Control Systems Handbook: Control System Advanced Methods}.
\newblock CRC Press, second edition, 2011.

\bibitem{Lew86}
F.~L. Lewis.
\newblock A survey of linear singular systems.
\newblock {\em Circuits Systems Signal Process.}, 5(1):3--36, 1986.

\bibitem{Mae98}
R.~M\"{a}rz.
\newblock Criteria for the trivial solution of differential algebraic equations
  with small nonlinearities to be asymptotically stable.
\newblock {\em J. Math. Anal. Appl.}, 225(2):587--607, 1998.

\bibitem{MMS16}
C.~Mehl, V.~Mehrmann, and P.~Sharma.
\newblock Stability radii for linear {H}amiltonian systems with dissipation
  under structure-preserving perturbations.
\newblock {\em SIAM J. Matrix Anal. Appl.}, 37(4):1625--1654, 2016.

\bibitem{MehMW18}
C.~Mehl, V.~Mehrmann, and M.~Wojtylak.
\newblock Linear algebra properties of dissipative {H}amiltonian descriptor
  systems.
\newblock {\em SIAM J. Matrix Anal. Appl.}, 39(3):1489--1519, 2018.

\bibitem{MehMW20}
C.~Mehl, V.~Mehrmann, and M.~Wojtylak.
\newblock Distance problems for dissipative hamiltonian systems and related
  matrix polynomials.
\newblock {\em Linear Algebra Appl.}, 2020.

\bibitem{Mi09thesis}
C.~K. Mikkelsen.
\newblock {\em Numerical methods for large {L}yapunov equations}.
\newblock ProQuest LLC, Ann Arbor, MI, 2009.
\newblock Thesis (Ph.D.)--Purdue University.

\bibitem{OtPa15}
G.~Ottaviani and R.~Paoletti.
\newblock A geometric perspective on the singular value decomposition.
\newblock {\em Rend. Istit. Mat. Univ. Trieste}, 47:107--125, 2015.

\bibitem{ReReVo15}
T.~Reis, O.~Rendel, and M.~Voigt.
\newblock The {K}alman-{Y}akubovich-{P}opov inequality for
  differential-algebraic systems.
\newblock {\em Linear Algebra Appl.}, 485:153--193, 2015.

\bibitem{So98}
E.~D. Sontag.
\newblock {\em Mathematical control theory}, volume~6 of {\em Texts in Applied
  Mathematics}.
\newblock Springer-Verlag, New York, second edition, 1998.

\bibitem{St75}
T.~Str\"{o}m.
\newblock On logarithmic norms.
\newblock {\em SIAM J. Numer. Anal.}, 12(5):741--753, 1975.

\bibitem{Sty02}
T.~Stykel.
\newblock {\em Analysis and numerical solution of generalized {L}yapunov
  equations}.
\newblock PhD thesis, Technische Universit{\"a}t, Berlin, Institut f{\"u}r
  Mathematik, 2002.

\bibitem{Sty02a}
T.~Stykel.
\newblock Stability and inertia theorems for generalized {L}yapunov equations.
\newblock {\em Linear Algebra Appl.}, 355:297--314, 2002.

\bibitem{VD79}
P.~Van~Dooren.
\newblock The computation of {K}ronecker's canonical form of a singular pencil.
\newblock {\em Linear Algebra Appl.}, 27:103--140, 1979.

\bibitem{vD81}
P.~{Van Dooren}.
\newblock The generalized eigenstructure problem in linear system theory.
\newblock {\em IEEE Trans. Automat. Control}, 26(1):111--129, 1981.

\bibitem{Vi09}
C.~Villani.
\newblock Hypocoercivity.
\newblock {\em Mem. Amer. Math. Soc.}, 202(950):iv+141, 2009.

\bibitem{Wil72a}
J.~C. Willems.
\newblock Dissipative dynamical systems. {I}. {G}eneral theory.
\newblock {\em Arch. Rational Mech. Anal.}, 45:321--351, 1972.

\bibitem{Wil72b}
J.~C. Willems.
\newblock Dissipative dynamical systems. {II}. {L}inear systems with quadratic
  supply rates.
\newblock {\em Arch. Rational Mech. Anal.}, 45:352--393, 1972.

\bibitem{Wo85}
W.~M. Wonham.
\newblock {\em Linear multivariable control}, volume~10 of {\em Applications of
  Mathematics (New York)}.
\newblock Springer-Verlag, New York, third edition, 1985.

\end{thebibliography}

%

\end{document}